\documentclass[11pt]{amsart}
\usepackage{amsmath,amssymb,latexsym,dsfont}
\usepackage[left=2.6cm,right=2.6cm,top=2.7cm,bottom=2.7cm]{geometry}

\usepackage{amsmath,amsfonts,amssymb,epsfig, textcomp}
\usepackage{graphicx}
\usepackage{hyperref, cleveref}
\usepackage{cases,listings}
\usepackage{color}

\newtheorem{theorem}{Theorem}[section]
\newtheorem{lemma}{Lemma}[section]

\newtheorem{proposition}{Proposition}[section]
\newtheorem{corollary}{Corollary}[section]
\newtheorem{remark}{Remark}[section]

\newtheorem{definition}{Definition}[section]
\numberwithin{equation}{section}

      \newcommand{\hu}{\hat u}

      \newcommand{\hlambda}{\hat \lambda}

      \newcommand{\tf}{\widetilde f}

      \newcommand{\tlambda}{\widetilde \lambda}

 \newcommand{\cP}{{\mathcal P}}

   \newcommand{\tW}{\widetilde W}

   \newcommand{\bpsi}{\bar \psi}
   
      \newcommand{\tD}{\widetilde D}

      \newcommand{\tchi}{\widetilde \chi}

 \newcommand{\cQ}{\mathcal Q}

   \newcommand{\ty}{\widetilde y}
 \newcommand{\tu}{\widetilde u}

 \newcommand{\bbeta}{\bar \eta}

      \newcommand{\cN}{{\mathcal N}}

      \newcommand{\cL}{{\mathcal L}}

      \newcommand{\N}{\mathbb{N}}

      \newcommand{\loc}{\operatorname{loc}}
      
      \newcommand{\eps}{\varepsilon}
      \newcommand{\mR}{\mathbb{R}}
      \newcommand{\mZ}{\mathbb{Z}}

      \newcommand{\mC}{\mathbb{C}}

      \newcommand{\dsp}{\displaystyle}

      \newcommand{\M}{{\mathcal M}}

      \makeatletter
      \def\@setcopyright{}
      \def\serieslogo@{}
      \makeatother
   \newcommand{\tr}{^\mathsf{T}}

\title[Decay for the KdV equations]{Decay for the nonlinear KdV equations at critical lengths}

\author{Hoai-Minh Nguyen}
\address[Hoai-Minh Nguyen]{Ecole Polytechnique F\'ed\'erale de Lausanne, EPFL, 
\newline \indent SB, CAMA, Station 8,  CH-1015 Lausanne, Switzerland.}
\email{hoai-minh.nguyen@epfl.ch}

\begin{document}

\maketitle

\begin{abstract} We consider the nonlinear Korteweg-de Vries (KdV) equation in a bounded interval equipped with the Dirichlet boundary condition and the Neumann boundary condition on the right.  It is known that there is a set of critical lengths for which 
the solutions of the linearized system conserve the $L^2$-norm if their initial data belong to a finite dimensional  subspace $\M$. 
In this paper, we show that all solutions of the nonlinear KdV system decay to 0 at least with the rate $1/ t^{1/2}$ when $\dim \M = 1$ or when $\dim \M$ is even and a specific  condition is satisfied,  provided that  their initial data is sufficiently small. Our analysis is inspired by  the power series expansion approach and involves the theory of quasi-periodic functions. As a consequence, we rediscover known  results which were previously established   for $\dim \M = 1$ or for the smallest critical length $L$ with $\dim \M = 2$  by a different approach using the center manifold theory, and obtain new results.  We also show that the decay rate is not slower than $\ln (t + 2) / t$ for all critical lengths. 
\end{abstract}

\noindent{\bf Key words:} KdV equations, critical lengths, decay of solutions, asymptotically stable, quasi-periodic functions.

\noindent{\bf AMS subject classification:} 35B40, 35C20, 35Q53, 93B05.

\tableofcontents

\section{Introduction}

\subsection{Introduction and statement of the main results}
We consider the nonlinear Korteweg-de Vries (KdV) equation in a bounded interval $(0, L)$ equipped with the Dirichlet boundary condition and the Neumann boundary condition on the right: 
\begin{equation}\label{KdV-NL}\left\{
\begin{array}{cl}
u_t (t, x) + u_x (t, x) + u_{xxx} (t, x) + u (t,x) u_x (t, x)  = 0 &  \mbox{ for } t \in (0, +\infty), \, x \in (0, L), \\[6pt]
u(t, x=0) = u(t, x=L) = u_x(t , x= L)=  0 & \mbox{ for } t \in (0, +\infty), \\[6pt]
u(t = 0, \cdot)  = u_0 & \mbox{ in } (0, L), 
\end{array}\right.
\end{equation}
where $u_0 \in L^2(0, L)$ is the initial data.  The KdV equation has been introduced by Boussinesq \cite{1877-Boussinesq} and Korteweg
and de Vries \cite{KdV} as a model for propagation of surface water waves along a
channel.
This equation also furnishes a very useful nonlinear approximation model
including a balance between a weak nonlinearity and weak dispersive effects and has been studied extensively, see e.g.~\cite{Whitham74, Miura76}. 

Regarding \eqref{KdV-NL}, Rosier \cite{Rosier97} introduced a  set of critical lengths $\cN$ defined by 
\begin{equation}\label{def-cN}
\cN : = \left\{ 2 \pi \sqrt{\frac{k^2 + kl + l^2}{3}}; \, k, l \in \N_*\right\}. 
\end{equation}
This set plays an important role in both the decay property of the solution $u$ of \eqref{KdV-NL} and the controllability property of the system associated with \eqref{KdV-NL} where $u_x(t, L)$ is a control instead of $0$. 

Let us briefly review the known results on the controllability of \eqref{KdV-NL}  where $u_x(t, L)$ is a control: 
\begin{equation}\label{KdV-NLC}\left\{
\begin{array}{cl}
u_t (t, x) + u_x (t, x) + u_{xxx} (t, x) + u (t, x) u_x (t, x)  = 0 &  \mbox{ for } t \in (0, +\infty), \, x \in (0, L), \\[6pt]
u(t, x=0) = u(t, x=L) =   0 & \mbox{ for } t \in (0, +\infty), \\[6pt]
u_x(\cdot , x= L) : \mbox{ is a control}, \\[6pt]
u(t = 0, \cdot)  = u_0 & \mbox{ in } (0, L). 
\end{array}\right.
\end{equation}
For initial and final datum in $L^2(0, L)$ and controls in $L^2(0, T)$,
Rosier~\cite{Rosier97} proved that  system \eqref{KdV-NLC}
is small-time locally controllable around 0  provided that the length $L$ is not critical,
i.e., $L \not \in \cN$.  To this end, he studied the controllability of the linearized system using the Hilbert
Uniqueness Method  and compactness-uniqueness arguments.  He also established
that when the length $L$ is critical, i.e., $L \in \cN$,  the linearized system is not controllable. More precisely, he showed
that there exists a non-trivial finite-dimensional subspace $\M$ ($= \M_L$) of $L^2(0, L)$  such that its orthogonal space in $L^2(0, L) $ is reachable from $0$ whereas  $\M$ is not. To tackle the control problem for the critical length $L \in \cN$, Coron and Cr\'epeau
introduced the power series expansion method \cite{CC04}. The idea is to take into account the effect of
 the nonlinear term $u u_x$  absent in  the linearized system. Using this
method, they showed \cite{CC04} (see also \cite[section 8.2]{Coron07}) that  system \eqref{KdV-NLC}   is
small-time locally controllable if $L = m 2 \pi$ for $m \in \N_*$ satisfying
\begin{equation}
\nexists (k, l) \in \N_* \times \N_* \mbox{ with } k^2 + kl + l^2 = 3 m^2 \mbox{ and }
k \neq l, 
\end{equation}
with  initial and final  datum in $L^2(0, L)$ and controls in $L^2(0, T)$. 
In this case,  $\dim \M = 1$ and $\M$ is spanned by $1 - \cos x$. Cerpa \cite{Cerpa07}
developed the analysis in \cite{CC04} to prove that  \eqref{KdV-NLC}  is locally
controllable at \emph{a finite time} in the case $\dim \M = 2$. This corresponds to
the
case where
\[
L = 2 \pi \sqrt{\frac{k^2 + kl + l^2}{3}}
\]
for some $k, \,  l \in \N_*$ with   $k>l$, and there is no $(m, n) \in \N_* \times \N_*$ with $m>n$
and $m^2 + mn + n^2 = k^2 + kl + l^2$. Later, Cr\'epeau and Cerpa \cite{CC09}
succeeded to extend  the ideas in \cite{Cerpa07} to obtain the local
controllability for all other critical lengths at {\it a finite time}. Recently, Coron, Koenig, and Nguyen \cite{CKN-20} prove that when 
$(2k + l) / 3 \not \in \N_*$, one cannot achieve the small time local controllability for initial datum in $H^3(0, L)$ and controls in $H^1$ (in time). We also establish the local controllability for finite time of \eqref{KdV-NLC} for some subclass of these pairs $(k, l)$ with initial datum in $H^3(0, L)$ and the controls in $H^1(0, T)$. This is surprising when compared with known results on internal controls for  system \eqref{KdV-NL}. It is known, see \cite{CPR15, PVZ02, Pazoto05}, that system \eqref{KdV-NL} is locally controllable using internal controls {\it whenever} the  control region contains  an {\it arbitrary} open subset of $(0, L)$.

We next discuss the decay property of \eqref{KdV-NL}. Multiplying the equation of $u$ (real) with $u$ and integrating by parts, one obtains 
\begin{equation}\label{key-identity}
\int_{0}^L |u(t, x)|^2 \, dx + \int_0^t |u_x(s, 0)|^2 \, ds  = \int_{0}^L |u(0, x)|^2 \, dx \mbox{ for all } t > 0. 
\end{equation}
As a consequence of \eqref{key-identity}, one has
\begin{equation}\label{key-identity-0}
\int_{0}^L |u(t, x)|^2 \, dx \le \int_{0}^L |u(0, x)|^2 \, dx \mbox{ for all } t > 0. 
\end{equation}
In the case $L \not \in \cN$,  Menzala, Vasconcellos, and Zuazua  \cite{PVZ02} proved that the solutions of \eqref{KdV-NL} with small initial datum in $L^2(0, L)$ decay exponentially to 0. Their analysis is based on the exponential decay of the linearized system for which it holds, see \cite[Proposition 3.3]{Rosier97}, 
\begin{equation}\label{key-identity-1}
\int_0^t |u_x(s, 0)|^2 \, ds \ge c_t \int_{0}^L |u(0, x)|^2 \mbox{ for all } t > 0. 
\end{equation}
When a local damping was added, they also obtained the global exponential stability using the multiplier technique, compactness arguments, and the unique continuation for the KdV equations. 
Related results on modified nonlinear KdV equations can be found in \cite{RZ06,LP07}.  It is known from the work of Rosier \cite{Rosier97} that for $u_0 \in \M$, the solution $u$ of the linearized system satisfies 
\begin{equation}\label{key-identity-2}
\int_0^t |u_x(s, 0)|^2  \, ds  = 0 \mbox{ for all } t > 0,
\end{equation}
which implies in particular that  \eqref{key-identity-1} does not hold for any $t > 0$.  The work of Menzala, Vasconcellos, and Zuazua naturally raises the question whether or not  the solutions of \eqref{KdV-NL} go to 0 as the time goes to infinity  (see \cite[Section 4]{PVZ02} and also \cite[Section 5]{Pazoto05}). Quite recently, progress has been made for this problem.  Concerning the decay property of \eqref{KdV-NL} for critical lengths,  when $\dim \M = 1$,  Chu, Coron, and Shang \cite{CCS15} showed that the solution $u(t, \cdot)$ goes to 0 as $t \to+ \infty$ for all small initial data in $L^2(0, L)$. Moreover, they showed that there exists a constant $C$ depending only on $L$ such that 
\begin{equation}\label{CCS15}
\| u(t, \cdot) \|_{L^2(0, L)} \le \frac{C}{\sqrt{t}} \mbox{ for } t > 0. 
\end{equation}
It is worth mentioning that  the set of $L \in \cN$ such that $\dim \M = 1$ is infinite \cite{CC04}.  When $k = 2$ and $l = 2$ (the smallest length for which $\dim \M  = 2$), Tang, Chu, Sang, and Coron \cite{TCSC18} also established the decay to 0 of the solutions
by establishing an estimate equivalent to \eqref{CCS15} (see \cite[(1.20) in Theorem 1.1]{TCSC18}).  The analysis in \cite{CCS15,TCSC18} is based on the center manifold theory in infinite dimensions, see e.g. \cite{HI11}, in particular the work \cite{VMW04}. To this end, the authors showed the existence and smoothness of a center manifold associated with \eqref{KdV-NL}, which have their own interests. 

\medskip 
In this paper, we show  that  all solutions of \eqref{KdV-NL} decay to 0 at least with  a rate $1/t^{1/2}$ provided their initial data in $L^2(0, L)$ is small enough when $\dim \M = 1$ or when  condition \eqref{main-assumption} below holds (this requires in particular that $\dim \M$ is even). Given 
a critical length $L$, condition~\eqref{main-assumption}  can be checked numerically, a scilab program is given  in the appendix (see \Cref{cor1} for a range of validation).   Our approach is inspired by the spirit of the power series expansion due to Coron and Cr\'epeau \cite{CC04} and involves the theory of quasi-periodic functions.  

Before stating our results, let us introduce some notations  associated with  the structure of $\M$, see e.g. \cite{Rosier97,CC04,Cerpa14}. Recall that, for each $L \in \cN$,  there exists exactly $n_L \in \N_*$ pairs $(k_m, l_m) \in \N_* \times \N_*$ ($1 \le m \le n_L$) such that $k_m \ge l_m$,  and 
\begin{equation}\label{def-L}
L = 2 \pi \sqrt{\frac{k_m^2 + k_m l_m + l_m^2}{3}}. 
\end{equation}
For $1 \le m \le n_L$, set 
\begin{equation}\label{def-pm}
p_m = p(k_m, l_m) = \frac{(2k_m + l_m)(k_m - l_m)(2 l_m + k_m)}{3 \sqrt{3}(k_m^2 + k_m l_m + l_m^2)^{3/2}}, 
\end{equation}
 and   denote 
\begin{equation}
\cP_L = \Big\{p_m \mbox{ given by } \eqref{def-pm}; 1 \le m \le n_L \Big\}. 
\end{equation}
For $L \in \cN$ and $1 \le m \le n_L$ with $p_m \neq 0$, let $\sigma_{j, m}$ ($1 \le j \le 3$) be the solutions of 
$$
\sigma^3  - 3 (k_m^2 + k_m l_m + l_m^2) \sigma + 2(2 k_m + l_m)(2 l_m + k_m) (k_m - l_m) = 0, 
$$
and set, with the convention $\sigma_{j+3, m} = \sigma_{j, m}$ for $j \ge 1$, 
\begin{equation}\label{def-sm}
s_m = s(k_m, l_m) : = \sum_{j=1}^3 \sigma_{j, m} (\sigma_{j+2, m}  - \sigma_{j+1, m} ) \left( e^{\frac{4 \pi i (k_m - l_m)}{ 3} }   e^{2 \pi i \sigma_{j, m}} +  e^{- 2 \pi i \sigma_{j, m}}\right). 
\end{equation}

We are ready to state the main result of the paper: 

\begin{theorem}\label{thm1} Let $L \in \cN$.  Assume that either $\dim \M = 1$ or
\begin{equation}\label{main-assumption}
p_m \neq 0 \quad \mbox{ and } \quad s_m \neq 0  \quad \mbox{ for all } 1 \le m \le n_L. 
\end{equation}
There exists $\eps_0 > 0$ depending only on $L$ such that for all (real) $u_0 \in L^2(0, L)$ with $\| u_0 \|_{L^2(0, L)} \le \eps_0$, the unique solution $u \in C \big([0, + \infty); L^2(0, L) \big) \cap L^2_{\loc} \big([0, + \infty); H^1(0, L) \big)$
of \eqref{KdV-NL} satisfies 
\begin{equation}\label{thm1-cl1}
\lim_{t \to 0} \| u(t, \cdot) \|_{L^2(0, L)}  = 0. 
\end{equation}
More precisely, there exists a constant $C$ depending only on $L$ such that, for  $t \ge C/ \| u_0 \|_{L^2(0, L)}^2$ and $\| u_0 \|_{L^2(0, L)} \le \eps_0$, it holds  
\begin{equation}\label{thm1-cl2}
\|u(t, \cdot)  \|_{L^2(0, L) } \le \frac{1}{2} \| u(0, \cdot) \|_{L^2(0, L)}. 
\end{equation}
As a consequence, we have  
\begin{equation}\label{thm1-cl3}
\| u(t, \cdot) \|_{L^2(0, L)} \le c /t^{1/2} \mbox{ for } t \ge 0, 
\end{equation}
for some positive constant $c$ depending only on $L$. 
\end{theorem}

\begin{remark} \rm Let $L \in \N$. Condition $p_m \neq 0$ for all $1 \le m \le n_L$ is equivalent to the fact that  $\dim \M$ is even, see e.g.  \cite{Cerpa14}.

\begin{remark}\rm  Note that $s_m$ is a antisymmetric function of $(\sigma_{1, m}, \sigma_{2, m}, \sigma_{3, m})$ and hence the condition \eqref{main-assumption} does not depend on the order of $(\sigma_{1, m}, \sigma_{2, m}, \sigma_{3, m})$. 
\end{remark}

\begin{remark} \rm Assume  \eqref{main-assumption}. Applying \Cref{thm1}, one derives from  \eqref{key-identity-0} that $0$ is (locally) asymptotically stable with respect to $L^2(0, L)$-norm  for system  \eqref{KdV-NL}. 
\end{remark}

\begin{remark} \rm Assume that \eqref{thm1-cl2} holds. By the regularity properties of the KdV equations, one derives that the same rate of decay holds for $t>1$ when $\| \cdot \|_{L^2(0, L) }$ is replaced by $\| \cdot \|_{H^m(0, L) }$ for $m \ge 1$. 
\end{remark}

\end{remark}

Condition \eqref{main-assumption} can be checked numerically. For example, using scilab (the program is given in the appendix), we can check $s_m \neq 0$ for all $(k_m, l_m) \in \N_*$ with $1 \le l_m < k_m < 2000$. As a consequence, we have 

\begin{corollary}\label{cor1} Let $L \in \cN$. Assume that either $\dim \M = 1$ or $1 \le k_m, l_m \le 1000$ for some $1 \le m \le n_L$. Then \eqref{thm1-cl3} holds if $p_m \neq 0$ for all $1 \le m \le n_L$. 
\end{corollary}

We thus rediscover  the decay results in \cite{CCS15,TCSC18} by a different approach and  obtain new results.

\begin{remark} \rm
Concerning \eqref{main-assumption}, we expect that  $s_m \neq 0$ holds for all $L \in \cN$ but  we are not able to show it. 
\end{remark}

The optimality of the decay rate $1/ t^{1/2}$ given in \eqref{thm1-cl3} is open. However, we can establish the following result for all critical lengths. 

\begin{proposition}\label{pro-opt} Let $L \in \cN$. There exists $c > 0$ such that for all $\eps > 0$, there exists $u_0 \in L^2(0, L)$ such that 
$$
\| u_0\|_{L^2(0, L)} \le \eps \quad \mbox{ and } \quad \| u(t, \cdot) \|_{L^2(0, L)} \ge  c \ln (t+2)  /t \mbox{ for some } t > 0.  
$$
\end{proposition}

It is natural  to ask if the decay  holds globally, i.e., without the assumption on the smallness of the initial data. In fact, this cannot hold even for non-critical lengths.  More precisely, Doronin and Natali \cite{DN14} showed that there exist (infinite) stationary states of \eqref{KdV-NL}  for any  length $L$,  which is critical or not.

\subsection{Ideas of the analysis and structure of the paper} The key of the analysis of \Cref{thm1} is to (observe and) establish the following fact (see \Cref{lemK}): Let $L \in \cN$. Under condition \eqref{main-assumption} or $\dim \M =1$,  there exist two constants $T_0>0$ and $C>0$ depending only on $L$ such that for $T \ge T_0$, one has, for all $u_0 \in L^2(0, L)$ with  $\| u_0\|_{L^2(0, L)}$ sufficiently small,  
\begin{equation}\label{decayK-I}
\| u(T, \cdot) \|_{L^2(0, L)} \le  \| u_0 \|_{L^2(0, L)} \Big(1 -  C \| u_0 \|_{L^2(0, L)}^2 \Big) \mbox{ for } T \ge T_0,  
\end{equation}
where $u$ is the unique solution of \eqref{KdV-NL}. 

To get an idea of how to prove \eqref{decayK-I}, let us consider the case $u_0 \in \M \setminus \{0 \}$, which is somehow the worst case.  The analysis is inspired by  the spirit  of the power expansion method   \cite{CC04}. Let  $\tu_1$ be  the unique solution of
\begin{equation} \label{hu1-Int}\left\{
\begin{array}{cl}
\tu_{1, t} (t, x) + \tu_{1, x} (t, x) + \tu_{1, xxx} (t, x)   = 0 &  \mbox{ for } t \in (0, +\infty), \, x \in (0, L), \\[6pt]
\tu_1(t, x=0) = \tu_1(t, x=L) = \tu_{1, x} (t , x= L)=  0 & \mbox{ for } t \in (0, +\infty), \\[6pt]
\tu_1(t = 0, \cdot)  = u_0/ \eps & \mbox{ in } (0, L), 
\end{array}\right.
\end{equation}
with $\eps = \| u_0\|_{L^2(0, L)} > 0$, and let $\tu_2$ be the unique solution of 
\begin{equation}\label{hu2-Int}\left\{
\begin{array}{cl}
\tu_{2, t} (t, x) + \tu_{2, x} (t, x) + \tu_{2, xxx} (t, x) + \tu_{1, x} (t, x) \tu_1 (t, x)   = 0 &  \mbox{ for } t \in (0, +\infty), \, x \in (0, L), \\[6pt]
\tu_2(t, x=0) = \tu_2(t, x=L) = \tu_{2, x} (t , x= L)=  0 & \mbox{ for } t \in (0, +\infty), \\[6pt]
\tu_2(t = 0, \cdot)  = 0  & \mbox{ in } (0, L).  
\end{array}\right.
\end{equation}
By considering the system of $\eps \tu_1 + \eps^2 \tu_2 - u$, we can prove that, for arbitrary $\tau > 0$,  
\begin{equation}\label{diff-I}
\| (\eps \tu_1 + \eps^2 \tu_2 - u)_x (\cdot, 0) \|_{L^2(0, \tau)} \le c_\tau \eps^3, 
\end{equation}
for some $c_\tau > 0$ depending only on $\tau$ and $L$, provided that  $\eps$ is sufficiently small. Since $\tu_{1}(t, \cdot) \in \M$ for all $t > 0$, one can then derive that 
$$
\tu_{1, x} (t, 0) = 0 \mbox{ for } t \ge 0.  
$$
Thus, if one can show that, for some $\tau_0 > 0$ and for some  $c_0> 0$ 
\begin{equation}\label{cond-hu2-I}
\|  \tu_{2, x}(\cdot, 0) \|_{L^2(0, \tau_0)} \ge c_0, 
\end{equation}
then from \eqref{diff-I} one has, for $\eps$ small enough, 
$$
\|  u_{x}(\cdot, 0) \|_{L^2(0, \tau_0)} \ge c_0 \eps^2. 
$$
This implies \eqref{decayK-I} with $T_0 = \tau_0$ by \eqref{key-identity}. 

To establish \eqref{cond-hu2-I}, we first construct a special solution $W$ of the system 
\begin{equation}\label{W-Int}\left\{
\begin{array}{cl}
W_t (t, x) + W_{x} (t, x) + W_{xxx} (t, x) + \tu_{1, x} (t, x) \tu_1 (t, x)   = 0 &  \mbox{ for } t \in (0, +\infty), \, x \in (0, L), \\[6pt]
W(t, x=0) = W(t, x=L) = W (t , x= L)=  0 & \mbox{ for } t \in (0, +\infty), 
\end{array}\right.
\end{equation}
via  a separation-of-variable process. Moreover, we can prove for such a solution $W$ that 
\begin{multline}\label{W_x-quasi-I}
\mbox{$W$ is bounded by $\| \tu_1(0, \cdot) \|_{L^2(0, L)}$ up to a  positive constant,}  \\
\mbox{and  $W_{x}(\cdot, 0)$ is a non-trivial quasi-periodic function.}
\end{multline}
The proof of this property is based on some useful observations on $p_m$ and the boundary conditions considered in \eqref{KdV-NL},  and involves some arithmetic arguments.  It is in the proof of the existence of $W$ and  the second fact of \eqref{W_x-quasi-I} that  assumption  \eqref{main-assumption} or $\dim \M = 1$ is required.  Note that, for all $\delta > 0$, there exists $T_\delta > 0$ such that it holds, for $\tau \ge T_\delta$, 
\begin{equation}\label{decay-I}
\| y_x(\cdot, 0) \|_{L^2(\tau, 2 \tau)} \le \delta \| y_0 \|_{L^2(0, L)}, 
\end{equation}
for all solution $y \in C \big([0, + \infty); L^2(0, L) \big) \cap L^2_{\loc} \big([0, + \infty); H^1(0, L) \big)$ of the system 
\begin{equation*}\left\{
\begin{array}{cl}
y_t (t, x) + y_x (t, x) + y_{xxx} (t, x)   = 0 &  \mbox{ for } t \in (0, +\infty), \, x \in (0, L), \\[6pt]
y(t, x=0) = y(t, x=L) = y_x(t , x= L)=  0 & \mbox{ for } t \in (0, +\infty). 
\end{array}\right.
\end{equation*}
Combining \eqref{W_x-quasi-I} and \eqref{decay-I}, we can derive \eqref{cond-hu2-I} after applying the theory of quasi-periodic functions, see e.g. \cite{Bohr47}. 

\medskip  The proof of \Cref{pro-opt}  is inspired by the approach which is used to prove  \Cref{thm1} and is mentioned above. 

\medskip 
The paper is organized as follows. The elements for the construction of $W$ are given in \Cref{sect-construction} and the elements for the
proof of \eqref{W_x-quasi-I} are given in \Cref{sect-quasi}. The proof of \Cref{thm1} is given in \Cref{sect-thm1} where 
\eqref{decayK-I} is formulated in \Cref{lemK}. The proof of \Cref{pro-opt} is given in \Cref{sect-opt}. 
In the appendix, we reproduce a proof of a technical result, which is obtained in \cite{CKN-20},  and provide the scilab code.

\section{Construction of auxiliary functions} \label{sect-construction}

Let us begin with recalling and introducing some useful notations motivated by the structure of $\M$, see e.g. \cite{Rosier97,CC04,Cerpa14}.  For $L \in \cN$ and for $1 \le m \le n_L$, denote 
\begin{equation}\label{def-etam}
\left\{\begin{array}{c}
\dsp \eta_{1, m} = - \frac{2 \pi i  (2 k_m + l_m) }{3 L },\\[6pt]  
\dsp \eta_{2, m} = \eta_{1, m} + \frac{2 \pi i }{L} k_m =  \frac{2 \pi i  (k_m - l_m) }{3 L }, \\[6pt] 
\dsp \eta_{3, m} =  \eta_{2, m} + \frac{2 \pi i }{L} l_m =    \frac{2 \pi i  (k_m + 2  l_m) }{3 L }.  
\end{array}\right.
\end{equation}
Set  
\begin{equation}\label{def-psi}
\left\{ \begin{array}{cl}
 \psi_m(x) =  \sum_{j=1}^3 (\eta_{j+1, m} - \eta_{j, m}) e^{\eta_{j+2, m} x}  &  \mbox{ for } x \in [0, L], \\[6pt] \Psi_m(t, x) =  e^{- i t p_m}  \psi_m(x)  & \mbox{ for } (t, x) \in \mR \times [0, L], 
\end{array} \right.
\end{equation}
(recall that $p_m$ is defined in \eqref{def-pm}).  It is clear from the definition of $\eta_{j, m}$ in \eqref{def-etam} that
\begin{equation}\label{pro-etam}
e^{\eta_{1, m} L} = e^{\eta_{2, m} L} = e^{\eta_{3, m} L}. 
\end{equation} 
This property of $\eta_{j,m}$  associated with $L$ is used several times in our analysis.

\begin{remark} \rm  One can check that $\eta_{j, m}$ are the solutions of the equation 
$$
\lambda^3 + \lambda - i p_m \lambda = 0. 
$$
This implies in particular that $p_{m_1} \neq p_{m_2}$ if $(k_{m_1}, l_{m_1}) \neq (k_{m_2}, l_{m_2})$ as observed in \cite{Cerpa07}. 
\end{remark}

It is known that $\Psi_m$ is a solution of the linearized KdV system; moreover, 
$$
\Psi_{m, x}(\cdot, 0) = 0,  
$$
i.e., 
\begin{equation}\label{KdV-Psi}\left\{
\begin{array}{cl}
\Psi_{m, t} (t, x) + \Psi_{m, x} (t, x) + \Psi_{m, xxx} (t, x)  = 0 &  \mbox{ for } t \in (0, +\infty), \, x \in (0, L), \\[6pt]
\Psi_m(t, 0) = \Psi_m(t, L) = \Psi_{m, x} (t , 0) =\Psi_{m, x}(t , L) =   0 & \mbox{ for } t \in (0, +\infty). 
\end{array}\right.
\end{equation}
These properties  of $\Psi_m$ can be easily checked. It is known that, see e.g. \cite{Cerpa14}, 
\begin{equation}\label{span-M}
\M  = \mbox{span} \Big\{ \big\{ \Re( \psi_m(x)); 1 \le m \le n_L \Big\} \cup \Big\{ \Im( \psi_m(x)); 1 \le m \le n_L \big\} \Big\}. 
\end{equation}
Here and in what follows, for a complex number $z$, we denote $\Re z$, $\Im z$, and $\bar z$  its real part, its imaginary part, and its conjugate, respectively.

In this section, we prepare elements to construct the function $W$ mentioned in the introduction. Assume that $u_0 \in \M \setminus \{0\}$ and let $\eps = \| u_0 \|_{L^2(0, L)}$. By \eqref{span-M}, there exists $(\alpha_m)_{m=1}^{n_L} \subset \mC$ such that 
\begin{equation}\label{tu1-***}
\frac{1}{\eps} u_{0} = \Re \left\{ \sum_{m=1}^{n_L}  \alpha_m  \Psi_m (0, x) \right\}. 
\end{equation}
The function $\tu_1$ defined by \eqref{hu1-Int} is  then given by 
$$
\tu_1(t, x) = \Re \left\{ \sum_{m=1}^{n_L}  \alpha_m  \Psi_m (t, x) \right\} = \Re \left\{ \sum_{m=1}^{n_L}  \alpha_m  e^{- i p_m t}  \psi_m (x) \right\}. 
$$
Using the fact, for an appropriate complex function $f$, 
$$
\Re f(t, x) \Re f_x(t, x) = \frac{1}{2} \Big( (\Re f(t, x) )^2 \Big)_x= \frac{1}{8} \Big(  \big( f(t, x)^2)_x  + \big( \bar f(t, x) ^2\big)_x   + 2 ( |f(t, x)|^2)_x \Big), 
$$
we derive from \eqref{def-psi} and \eqref{tu1-***} that 
\begin{align}\label{motivation-sect-construction}
\tu_{1, x} (t, x) \tu_1(t, x) = & \frac{1}{8} \sum_{m_1=1}^{n_L}   \sum_{m_2=1}^{n_L}  \Big( \alpha_{m_1} \alpha_{m_2}  e^{-i (p_{m_1} + p_{m_2}) t} \psi_{m_1} (x) \psi_{m_2} (x) \Big)_x \\[6pt]
& +  \frac{1}{8} \sum_{m_1=1}^{n_L}  \sum_{m_2=1}^{n_L}   \Big( \overline{ \alpha_{m_1}  \alpha_{m_2}  e^{-i (p_{m_1} + p_{m_2}) t}  \psi_{m_1} (x)  \psi_{m_2} (x)} \Big)_x  \nonumber \\[6pt]
& +  \frac{1}{4} \sum_{m_1=1}^{n_L} \sum_{m_2=1}^{n_L} \Big( \alpha_{m_1} \bar \alpha_{m_2} e^{-i (p_{m_1} - p_{m_2}) t}   \psi_{m_1} (x) \bar \psi_{m_2} (x) \Big)_x.  \nonumber
\end{align}
Motivated by \eqref{motivation-sect-construction}, in this section, we construct solutions of system \eqref{pro1-sys-1}-\eqref{pro1-sys-2} and system \eqref{pro1-sys-1-Co}-\eqref{pro1-sys-2-Co} below. 

\medskip 
We begin with the following simple result whose proof is omitted. 

\begin{lemma}\label{lem-der} Let $L \in \cN$ and  $1 \le m_1, m_2 \le n_L$. We have, in $[0, L]$,  
\begin{multline}\label{lem-der-cl1}
\Big(\psi_{m_1} \psi_{m_2} \Big)'(x) \\[6pt]
= \sum_{j=1}^3 \sum_{k=1}^3 (\eta_{j+1, m_1} - \eta_{j, m_1})(\eta_{k+1, m_2} - \eta_{k, m_2}) (\eta_{j+2, m_1} + \eta_{k+2, m_2}) e^{(\eta_{j+2, m_1} + \eta_{k+2, m_2}) x}, 
\end{multline}
and 
\begin{multline}\label{lem-der-cl2}
\Big(\psi_{m_1} \bar \psi_{m_2} \Big)'(x) \\[6pt]
= \sum_{j=1}^3 \sum_{k=1}^3 (\eta_{j+1, m_1} - \eta_{j, m_1})(\bbeta_{k+1, m_2} - \bbeta_{k, m_2}) (\eta_{j+2, m_1} + \bbeta_{k+2, m_2}) e^{(\eta_{j+2, m_1} + \bbeta_{k+2, m_2}) x} . 
\end{multline}
\end{lemma}

We next introduce

\begin{definition}\label{def1}
For $z \in \mC$, let  $\lambda_j = \lambda_{j} (z)$ $(1 \le j \le 3)$ be the roots of the equation 
\begin{equation}\label{def-lambda}
\lambda^3 + \lambda -  i z  = 0,   
\end{equation}
and set 
\begin{equation}\label{def-Q}
Q (z) = \left(\begin{array}{ccc}
1 & 1 & 1 \\[6pt]
e^{\lambda_1 L } & e^{\lambda_2 L} & e^{\lambda_3 L} \\[6pt]
\lambda_1 e^{\lambda_1 L} & \lambda_2  e^{\lambda_2 L} & \lambda_3 e^{\lambda_3 L}
\end{array}\right).   
\end{equation}
\end{definition}

\begin{remark}\rm
Some comments on the definition of $Q$ are in order. 
The matrix $Q$ is antisymmetric  with respect
to $\lambda_j$ ($j=1, 2, 3$), and its definitions depend on a choice of  the order of  $(\lambda_1, \lambda_2, \lambda_3)$. Nevertheless, we later consider either the equation $\det Q = 0$ or a quantity  depending on $Q$ in such a way that  the order of  $(\lambda_1, \lambda_2, \lambda_3)$ does not matter. The definition of $Q$ is only considered in these contexts. 
\end{remark}

\begin{remark} \rm The definition of $\lambda_j(z)$ in \Cref{def1} is slightly  different from the one given in \cite{CKN-20} where $i z$ is used instead of $-iz$ in \eqref{def-lambda}. 
\end{remark}

\begin{remark} \rm \label{rem-pm-lambda} It is known that if $z \in \cP_L$ for some $L \in \cN$,  then 
$$
\lambda_j = \eta_{j, m} \mbox{ for some } 1 \le m \le n_L. 
$$
Hence, by \eqref{pro-etam}, 
$$
e^{\lambda_1 L  } = e^{\lambda_2 L  } = e^{\lambda_3 L }. 
$$
\end{remark}

\begin{remark}\label{rem-lambda} \rm Note that \eqref{def-lambda} has simple roots for $z \neq \pm 2/(3 \sqrt{3})$. Thus, a general solution of the equation 
$$
y'''(x) + y'(x) - i z y (x) = 0 \mbox{ in } [0, L], 
$$
is of the form $\sum_{j=1}^3 a_j e^{\lambda_j(z) x}$ when $z \neq \pm 2/ (3 \sqrt{3})$.  
For $z = \pm 2/(3 \sqrt{3})$, equation \eqref{def-lambda}  has three roots 
$$
\lambda_1 = \mp 2i / \sqrt{3} \quad \mbox{ and } \quad \lambda_2 = \lambda_3 = \pm i/ \sqrt{3}. 
$$
\end{remark}

We now recall a useful property of solutions of the equation $\det Q = 0$ which is established in \cite{CKN-20} 
(a consequence of  \cite[Remark 2.7]{CKN-20}).

\begin{lemma}\label{lem-Q}  Let $z \in \mR$.  Then $\det Q(z) = 0$ if and only if either $z = \pm 2/ \sqrt{3}$ or ($L \in \cN$ and $z \in \cP_L$). Moreover, 
$$
\big\{ \pm 2/ \sqrt{3} \big\}  \cap \cP_{L} = \emptyset \mbox{ for all } L \in \cN. 
$$
\end{lemma}

The proof of \Cref{lem-Q} is reproduced  in the appendix for the convenience of the reader.

\medskip 

Let $L \in \cN$ and  $1 \le m_1, m_2 \le n_L$. As mentioned above, we are interested in constructing  a solution of the system 
 \begin{equation}\label{pro1-sys-1}
- i (p_{m_1} + p_{m_2}) \varphi_{m_1, m_2}(x) + \varphi_{m_1, m_2}' (x) + \varphi_{m_1, m_2}''' (x) + \Big(\psi_{m_1} \psi_{m_2} \Big)'(x)= 0  \mbox{ in } (0, L), 
\end{equation}
and 
 \begin{equation}\label{pro1-sys-2}
\varphi_{m_1, m_2}(0) = \varphi_{m_1, m_2}(L) = \varphi_{m_1, m_2}'(L)=  0. 
\end{equation}
 
\medskip 
We have
 
\begin{proposition}\label{pro1} Let $L \in \cN$ and  $1 \le m_1, m_2 \le n_L$. Let $\lambda_j = \lambda_j(p_{m_1} + p_{m_2})$ and $\cQ = Q(i p_{m_1} + i p_{m_2})$ where $\lambda_j$ and $Q$ are defined by \eqref{def-lambda} and \eqref{def-Q}.  When $p_{m_1} \neq 0$ and $p_{m_2} \neq 0$, set
\begin{equation}\label{pro1-D}
D  = D_{m_1, m_2}= \sum_{j=1}^3 \sum_{k=1}^3 \frac{(\eta_{j+1, m_1} - \eta_{j, m_1})(\eta_{k+1, m_2} - \eta_{k, m_2}) }{3 \eta_{j+2, m_1}   \eta_{k+2, m_2}}, 
\end{equation}
and 
\begin{equation}\label{pro1-fm1m2}
\chi_{m_1, m_2}(x) = - \sum_{j=1}^3 \sum_{k=1}^3 \frac{(\eta_{j+1, m_1} - \eta_{j, m_1})(\eta_{k+1, m_2} - \eta_{k, m_2}) }{3 \eta_{j+2, m_1}   \eta_{k+2, m_2}} e^{(\eta_{j+2, m_1} + \eta_{k+2, m_2}) x} \mbox{ in } [0, L]. 
\end{equation}
We have 
\begin{enumerate}

\item[1)] Assume that $p_{m_1} \neq 0$,  $p_{m_2} \neq 0$,  and $p_{m_1} + p_{m_2}  \not \in \cP_L  \cup \big\{ 2 / (3 \sqrt{3}) \big\}$.  
The unique solution of system \eqref{pro1-sys-1}-\eqref{pro1-sys-2} is given by 
\begin{equation}\label{pro1-wm1m2}
\varphi_{m_1, m_2} (x) =  \chi_{m_1, m_2}(x) + \sum_{j=1}^3 a_j e^{\lambda_j x}, 
\end{equation}
where   $(a_1, a_2, a_3) $ is uniquely determined via \eqref{pro1-sys-2}, i.e., 
\begin{equation}\label{pro1-aj}
\cQ (a_1, a_2, a_3)\tr =    D(1, e^{(\eta_{1, m_1} + \eta_{1, m_2}) L }, 0)\tr.  
\end{equation}

\item[2)] Assume that $p_{m_1} \neq 0$,  $p_{m_2} \neq 0$,  and $p_{m_1} + p_{m_2}   \in \cP_L$.  A  solution of system \eqref{pro1-sys-1}-\eqref{pro1-sys-2} 
is given by \eqref{pro1-wm1m2} where $(a_1, a_2, a_3)$ satisfies
\begin{equation}\label{pro1-aj-2}
a_1 + a_2 + a_3 = D \quad \mbox{ and } \quad  \lambda_1 a_1 + \lambda_2 a_2 + \lambda_3 a_3 = 0. 
\end{equation}

\item[3)]  Assume that $p_{m_1} \neq 0$,  $p_{m_2} \neq 0$,  and $p_{m_1} + p_{m_2}  =  2/ (3 \sqrt{3})$.  Consider the convention  
\begin{equation}\label{lem2-lambda}
\lambda_1 = - 2 i / \sqrt{3} \quad \mbox{ and } \quad \lambda_2  = \lambda_3 = i / \sqrt{3}. 
\end{equation}
System \eqref{pro1-sys-1}-\eqref{pro1-sys-2} has a unique solution 
given by 
\begin{equation}\label{def-wmm-*}
\varphi_{m_1, m_2} (x) =  \chi_{m_1, m_2}(x) + 
   a_1 e^{\lambda_1 x} + (a_2 + a_3 x) e^{\lambda_2 x},  
\end{equation}
where $(a_1, a_2, a_3)$ is uniquely determined via \eqref{pro1-sys-2}, i.e., 
\begin{equation}\label{pro1-aj-3}
\cQ_1 (a_1, a_2, a_3)\tr =    D(1, e^{(\eta_{1, m_1} + \eta_{1, m_2}) L }, 0)\tr, 
\end{equation}
where 
\begin{equation}\label{def-Qm}
\cQ_1  = \left(\begin{array}{ccc}
1 & 1 & 0 \\[6pt]
e^{\lambda_1 L } & e^{\lambda_2 L} & L e^{\lambda_2 L} \\[6pt]
\lambda_1 e^{\lambda_1 L} & \lambda_2  e^{\lambda_2 L} &(\lambda_2 L + 1)  e^{\lambda_2 L}
\end{array}\right). 
\end{equation}

\item[4)] Assume that $p_{m_1}  =  p_{m_2} = 0$ and thus $m_1 = m_2 = m$.  A solution of system \eqref{pro1-sys-1}-\eqref{pro1-sys-2} is  
\begin{equation}\label{def-wmm-0}
\varphi_{m, m} (x) =   4 \left( L \sin x + \frac{1}{6} -  x \sin x  - \frac{1}{6} \cos (2x) \right). 
\end{equation}
\end{enumerate}
\end{proposition}

\begin{proof} We proceed with the proof of 1), 2), 3), and 4) in Steps 1, 2, 3, and  4 below, respectively.

\medskip 
\noindent{\it Step 1}: Proof of 1). Since $\eta  = \eta_{j, m}$ $(1 \le j \le 3)$ is a root of  the equation 
\begin{equation*}
\eta^3 + \eta - ip_m = 0, 
\end{equation*}
it follows that 
\begin{equation*} 
\eta_{j, m_1} \neq -  \eta_{k, m_2}
\end{equation*}
(since otherwise $p_{m_1} = - p_{m_2}$ which is impossible), and  
\begin{equation*}
(\eta_{j, m_1} + \eta_{k, m_2})^3 + (\eta_{j, m_1} + \eta_{k, m_2}) - i (p_{m_1} + p_{m_2}) = 3 \eta_{j, m_1} \eta_{k, m_2} (\eta_{j,m_1} + \eta_{k, m_2}). 
\end{equation*}
Since $p_{m_1} \neq 0$ and  $p_{m_2} \neq 0$, we derive from \Cref{lem-der} that $\chi_{m_1, m_2}$ is a solution of \eqref{pro1-sys-1}. Since a general solution of the equation $\xi''' + \xi' = i (p_{m_1} + p_{m_2}) \xi$ is of the form $\sum_{j=1}^3 a_j e^{\lambda_j x}$ by \Cref{rem-lambda},  
it follows that 
\begin{equation}\label{pro1-gen-sol}
\mbox{ a general solution of \eqref{pro1-sys-1} is of the form  $\chi_{m_1, m_2} (x) + \sum_{j=1}^3 a_j e^{\lambda_j x}$}. 
\end{equation}
We have 
\begin{equation}\label{pro1-pro-chi}
- \chi_{m_1, m_2} (0) = D, \quad - \chi_{m_1, m_2} (L) \mathop{=}^{\eqref{pro-etam}} D e^{(\eta_{1, m_1} + \eta_{1, m_2}) L}, \quad \mbox{ and } \quad - \chi_{m_1, m_2, x} (L) \mathop{=}^{\eqref{pro-etam}} 0. 
\end{equation}
It follows that a function of the form $\chi_{m_1, m_2} (x) + \sum_{j=1}^3 a_j e^{\lambda_j x}$ satisfies \eqref{pro1-sys-2} if and only if 
\begin{equation*}
\sum_{j=1}^3 a_j = D, \quad \sum_{j=1}^3 a_j e^{\lambda_j L} = D e^{(\eta_{1, m_1} + \eta_{1, m_2}) L}, \quad   \sum_{j=1}^3 a_j  \lambda_j e^{\lambda_j L} = 0, 
\end{equation*}
which is equivalent to \eqref{pro1-aj}. 
Since $p_{m_1} + p_{m_2} \not \in \cP_{L} \cup \big\{2 /( 3 \sqrt{3}) \big\}$ and $p_{m_1} + p_{m_2} > 0$, it follows from \Cref{lem-Q} that  $\det \cQ \neq 0$. Therefore, one obtains 1). 

\medskip 
\noindent{\it Step 2:} Proof of 2). A solution of \eqref{pro1-sys-1} is  of the form $\chi_{m_1, m_2} (x) + \sum_{j=1}^3 a_j e^{\lambda_j x}$. 
This function satisfies \eqref{pro1-sys-2} if and only if, by \Cref{rem-pm-lambda} (recall that $p_{m_1} + p_{m_2} \in \cP_L$),   
\begin{equation*}
\sum_{j=1}^3 a_j = D, \quad e^{\lambda_1 L} \sum_{j=1}^3 a_j  \mathop{=}^{\eqref{pro-etam}} D e^{(\eta_{1, m_1} + \eta_{1, m_2}) L}, \quad   \sum_{j=1}^3 a_j  \lambda_j  \mathop{=}^{\eqref{pro-etam}} 0. 
\end{equation*}
This system has a solution if 
\begin{equation}\label{Step2-key}
e^{\lambda_1 L} = e^{(\eta_{1, m_1} + \eta_{1, m_2}) L}, 
\end{equation}
and a solution is given 
by  \eqref{pro1-wm1m2} where $(a_1, a_2, a_3)$ satisfies \eqref{pro1-aj-2}.

It remains to prove \eqref{Step2-key}. Assume,  for some $p_{m_3} \in \cP_L$, that 
\begin{equation}\label{pro1-S2-p}
p_{m_1} + p_{m_2} = p_{m_3}. 
\end{equation}
To establish \eqref{Step2-key},  it suffices to prove that, by \eqref{pro-etam} and \Cref{rem-pm-lambda},  
$$
e^{( \eta_{2, m_1} + \eta_{2, m_2} ) L } = e^{\eta_{2, m_3} L } 
$$
which is equivalent to the fact, by \eqref{def-etam},  
\begin{equation}\label{pro1-S2-mod}
\frac{k_{m_3} - l_{m_3}}{3} - \frac{k_{m_1} - l_{m_1}}{3} - \frac{k_{m_2} - l_{m_2}}{3} \in \mZ. 
\end{equation}
From \eqref{pro1-S2-p} and the definition of $p_m$ in \eqref{def-pm}, we have 
\begin{multline}\label{pro1-S2-p1}
(k_{m_3} - l_{m_3}) (2k_{m_3} + l_{m_3}) (2 l_{m_3} + k_{m_3}) \\[6pt]
= (k_{m_1} - l_{m_1}) (2k_{m_1} + l_{m_1}) (2 l_{m_1} + k_{m_1}) + (k_{m_2} - l_{m_2}) (2k_{m_2} + l_{m_2}) (2 l_{m_2} + k_{m_2}). 
\end{multline}
Since 
$$
(k_{m_j} - l_{m_j}) (2k_{m_j} + l_{m_j}) (2 l_{m_j} + k_{m_j}) = l_{m_j}  - k_{m_j} \mod 3, 
$$
It follows from \eqref{pro1-S2-p1} that 
$$
k_{m_3} - l_{m_3} = k_{m_1} - l_{m_1} + \big( k_{m_2} - l_{m_2} \big)  \mod 3, 
$$
which yields \eqref{pro1-S2-mod}. The proof of Step 2 is complete. 

\medskip 
\noindent{\it Step 3:} Proof of 3). A solution of \eqref{pro1-sys-1} is of the form $\chi_{m_1, m_2}(x) +  a_1 e^{\lambda_1 x} + (a_2 + a_3 x) e^{\lambda_2 x}$. This function satisfies \eqref{pro1-sys-2} if and only if, by \eqref{pro1-pro-chi},  
$$
a_1 + a_2 = D, \quad a_1 e^{\lambda_1 L } + a_2 e^{\lambda_2 L } + a_3 L e^{\lambda_2 L} = De^{(\eta_{1, m_1} + \eta_{1, m_2} ) L}, 
$$
and 
$$
a_1 \lambda_1 e^{\lambda_1 L } + a_2 \lambda_2 e^{\lambda_2 L } + a_3 (\lambda_2 L + 1) e^{\lambda_2 L} = 0, 
$$
which is equivalent to \eqref{pro1-aj-3}.

Hence, it suffices to prove that $\cQ_1$ is invertible. Replacing the third row  of $\cQ_1$ by itself minus $\lambda_2$ times the second row, we obtain 
\begin{equation}
\cQ_2 = \left(\begin{array}{ccc}
1 & 1 & 0 \\[6pt]
e^{\lambda_1 L } & e^{\lambda_2 L} & L e^{\lambda_2 L} \\[6pt]
(\lambda_1 - \lambda_2) e^{\lambda_1 L} & 0 &e^{\lambda_2 L}
\end{array}\right). 
\end{equation}
We have 
$$
\det \cQ_2  = e^{2 \lambda_2 L} - \big(1 - L (\lambda_1 - \lambda_2) \big) e^{(\lambda_1 + \lambda_2 ) L }. 
$$
Using \eqref{lem2-lambda}, we derive that $\det \cQ_2 = 0$ if and only if 
$$
e^{3 \lambda_2 L } = 1 + 3 \lambda_2 L. 
$$
Since the equation $e^{ix } = 1 + i x$ has only one solution $x = 0$ in the real line, one derives that $\det \cQ_2 \neq 0$. Therefore, $\cQ_1$ is invertible.  The proof of Step 3) is complete.

\medskip 
\noindent{\it Step 4:} Proof of 4). Since $ p_{m} =0$, it follows that $k_m = l_m$, and $L = 2 \pi k_m$. One then has 
$$
\eta_{1, m} =  -i, \quad \eta_{2, m} = 0, \quad \eta_{3, m} = i. 
$$
It follows from the definition of $\psi_m$ in  \eqref{def-psi} that 
\begin{equation}\label{psi-complex}
\psi_m(x) = 2 i  (\cos x - 1).
\end{equation}
This implies 
$$
\big( \psi_m^2(x) \big)_x =   8 (\cos x - 1) \sin x. 
$$
A straightforward computation gives the conclusion.

\medskip 
The proof of \Cref{pro1} is complete. 
\end{proof}

\begin{remark} \rm In the case, $p_{m_1} = 0$ and $p_{m_2} \neq 0$, one cannot construct a solution of \eqref{pro1-sys-1}-\eqref{pro1-sys-2} in general.  In fact, one can check that 
\begin{multline}
\chi_{m_1, m_2} (x) =  - \sum_{j=1, 2} \sum_{k=1}^3 \frac{(\eta_{j+1, m_1} - \eta_{j, m_1})(\eta_{k+1, m_2} - \eta_{k, m_2}) }{3 \eta_{j+2, m_1}   \eta_{k+2, m_2}} e^{(\eta_{j+2, m_1} + \eta_{k+2, m_2}) x}  \\[6pt]
-  \sum_{k=1}^3 \frac{(\eta_{1, m_1} - \eta_{3, m_1})(\eta_{k+1, m_2} - \eta_{k, m_2})   \eta_{k+2, m_2} }{3 {\eta_{k+2, m_2}}^2 + 1} x e^{ \eta_{k+2, m_2} x}
\end{multline}
is a solution of \eqref{pro1-sys-1}. However, 
$$
\chi_{m_1, m_2} (0) \neq e^{-\eta_{1, m_2} L} \chi_{m_1, m_2} (L)
$$
since, in general,  
$$
 \sum_{k=1}^3 \frac{(\eta_{k+1, m_2} - \eta_{k, m_2})   \eta_{k+2, m_2} }{3 {\eta_{k+2, m_2}}^2 + 1}  \neq 0. 
$$
Hence one cannot find $(a_1, a_2, a_3) \in \mC^3$ such that the function $
\chi_{m_1, m_2}(x) + \sum_{j=1}^3 a_j e^{\lambda_j x}$, with $\lambda_j  = \lambda_j (p_{m_2})$, 
verifies \eqref{pro1-sys-2}. 
\end{remark}

Let $L \in \cN$ and  $1 \le m_1, m_2 \le n_L$. We are next interested in constructing  a solution of the system 
 \begin{equation}\label{pro1-sys-1-Co}
- i (p_{m_1} - p_{m_2}) \phi_{m_1, m_2}(x) + \phi_{m_1, m_2}'(x) + \phi_{m_1, m_2}''' (x) + \Big(\psi_{m_1} \bpsi_{m_2} \Big)'(x)= 0  \mbox{ in } (0, L), 
\end{equation}
and 
 \begin{equation}\label{pro1-sys-2-Co}
\phi_{m_1, m_2}(0) = \phi_{m_1, m_2}(L) = \phi_{m_1, m_2}'(L)=  0. 
\end{equation}

We have 

\begin{proposition}\label{pro1-Co} Let $L \in \cN$ and  $1 \le m_1, m_2 \le n_L$. Let $\tlambda_j = \lambda_j(p_{m_1} - p_{m_2})$ 
 and $\widetilde \cQ = Q(i p_{m_1} - i p_{m_2})$ where $\lambda_j$ and $Q$ are defined by \eqref{def-lambda} and \eqref{def-Q}.  When $p_{m_1} \neq 0$ and $p_{m_2} \neq 0$, set
\begin{equation}\label{pro1-D-Co}
\tD = \tD_{m_1, m_2} = \sum_{j=1}^3 \sum_{k=1}^3 \frac{(\eta_{j+1, m_1} - \eta_{j, m_1})(\bbeta_{k+1, m_2} - \bbeta_{k, m_2}) }{3 \eta_{j+2, m_1}   \bbeta_{k+2, m_2}}
\end{equation}
and 
\begin{equation}\label{pro1-chi-Co}
\tchi_{m_1, m_2}(x) = - \sum_{j=1}^3 \sum_{k=1}^3 \frac{(\eta_{j+1, m_1} - \eta_{j, m_1})(\bbeta_{k+1, m_2} - \bbeta_{k, m_2}) }{3 \eta_{j+2, m_1}   \bbeta_{k+2, m_2}} e^{(\eta_{j+2, m_1} + \bbeta_{k+2, m_2}) x} \mbox{ in } [0, L]. 
\end{equation}
We have 
\begin{enumerate}

\item[1)] Assume that $p_{m_1} \neq 0$,  $p_{m_2} \neq 0$,  $p_{m_1} \neq p_{m_2}$, and $p_{m_1} - p_{m_2}  \not \in \cP_L $.  
The unique solution of system \eqref{pro1-sys-1-Co}-\eqref{pro1-sys-2-Co} is given by 
\begin{equation}\label{pro1-wm1m2-Co}
\phi_{m_1, m_2} (x) =  \tchi_{m_1, m_2}(x) 
 + \sum_{j=1}^3 a_j e^{\tlambda_j x}, 
\end{equation}
where   $(a_1, a_2, a_3)$ is uniquely determined via \eqref{pro1-sys-2-Co}, i.e., 
\begin{equation}\label{def-aj-Co}
\widetilde \cQ (a_1, a_2, a_3)\tr =    \tD(1, e^{(\eta_{1, m_1} + \bbeta_{1, m_2}) L }, 0)\tr.  
\end{equation}

\item[2)] Assume that $p_{m_1} \neq 0$,  $p_{m_2} \neq 0$, $p_{m_1} \neq p_{m_2}$, and $p_{m_1} - p_{m_2}   \in \cP_L$.  A  solution of system \eqref{pro1-sys-1-Co}-\eqref{pro1-sys-2-Co} 
is given by \eqref{pro1-wm1m2-Co} where $(a_1, a_2, a_3)$ satisfies  
\begin{equation}
a_1 + a_2 + a_3 = \tD \quad \mbox{ and } \tlambda_1 a_1 + \tlambda_2 a_2 + \tlambda_3 a_3 = 0. 
\end{equation}

\item[3)] Assume that $p_{m_1}  =  p_{m_2} \neq 0$ and thus $m_1 = m_2 = m$.  System \eqref{pro1-sys-1-Co}-\eqref{pro1-sys-2-Co} has a unique solution 
\begin{multline*}
\phi_{m, m} (x)  =   - \sum_{j=1}^3 \sum_{k=1}^3 \frac{(\eta_{j+1, m} - \eta_{j, m})(\bbeta_{k+1, m} - \bbeta_{k, m}) }{3 \eta_{j+2, m}   \bbeta_{k+2, m}} e^{(\eta_{j+2, m} + \bbeta_{k+2, m}) x} \\[6pt]
 +  \sum_{j=1}^3 \sum_{k=1}^3 \frac{(\eta_{j+1, m} - \eta_{j, m})(\bbeta_{k+1, m} - \bbeta_{k, m}) }{3 \eta_{j+2, m}   \bbeta_{k+2, m}}. 
\end{multline*}
\end{enumerate}

\item[4)] Assume that $p_{m_1}  =  p_{m_2} = 0$ and thus $m_1 = m_2 = m$.  A solution of  system \eqref{pro1-sys-1-Co}-\eqref{pro1-sys-2-Co} is  
\begin{equation}\label{def-wmm-0-Co}
\phi_{m, m} (x) =   - 4 \left( L \sin x + \frac{1}{6} -  x \sin x  - \frac{1}{6} \cos (2x) \right). 
\end{equation}

\end{proposition}

\begin{proof}  We proceed with the proof of 1), 2), 3), and  4) in Steps 1, 2, 3, and 4 below, respectively.

\medskip 
\noindent{\it Step 1}: Proof of 1). The proof is similar tos Step 1 in the proof of \Cref{pro1}. One just notes that 
\begin{equation*}
(\eta_{j, m_1} + \bbeta_{k, m_2})^3 + (\eta_{j, m_1} + \bbeta_{k, m_2}) - i (p_{m_1} - p_{m_2}) = 3 \eta_{j, m_1} \bbeta_{k, m_2} (\eta_{j, m_1} + \bbeta_{k, m_2}),  
\end{equation*}
and 
$$
\eta_{j, m_1} + \bbeta_{k, m_2} \neq 0
$$
since $p_{m_1} \neq p_{m_2}$. 

\medskip 
\noindent{\it Step 2:} Proof of 2). The proof is almost the same as  Step 2 in the proof of \Cref{pro1}. The details are omitted.

\medskip 
\noindent{\it Step 3:} Proof of 3). One can check that $\phi_{m, m}$ is a solution of \eqref{pro1-sys-1-Co}-\eqref{pro1-sys-2-Co}. The uniqueness follows from the fact that equation \eqref{def-lambda} has simple roots for $z=0$.

\medskip 
\noindent{\it Step 4:} Proof of 4). The conclusion is from 4) of \Cref{pro1} by noting that 
$$
|\psi_m(x)|^2 \mathop{=}^{\eqref{psi-complex}} - \psi_m(x)^2 \mbox{ if } p_m =0. 
$$

\medskip 
The proof is complete. 
\end{proof}

\section{Properties of auxiliary functions} \label{sect-quasi}

The main goal of this section is to establish,  for $L \in \cN$ and $1 \le m \le n_L$ with $p_m \neq 0$, that 
\begin{equation}\label{der-varphi-0}
\varphi_{m, m}'(0) \neq 0
\end{equation}
provided \eqref{main-assumption} holds (see \Cref{pro-chi}) where $\varphi_{m, m}$ is determined in \Cref{pro1}.  We begin with 

\begin{lemma} \label{lemE} Let $L \in \cN$ and  $1 \le m \le n_L$ with  $p_m \neq 0$. 
Set 
\begin{equation}
E_m  : = \sum_{j=1}^3 \frac{\eta_{j+1, m} - \eta_{j, m}}{\eta_{j+2, m}}.
\end{equation}
We have 
\begin{equation}\label{lemE-cl1}
D_{m, m} = - \chi_{m, m}(0) = \frac{1}{3} E_m^2, 
\end{equation}
and 
\begin{equation}\label{lemE-cl2}
E_m
=  - \frac{2 7 k_m l_m (k_m+l_m)}{(k_m+2l_m) (2 k_m + l_m) (k_m-l_m)} \neq 0. 
\end{equation}
\end{lemma}

\begin{proof}   It is clear to see from \eqref{pro1-fm1m2} that 
$$
D_{m, m} = - \chi_{m, m}(0) = \frac{1}{3} E_m^2. 
$$

With the notation $\gamma_{j, m} = L \eta_{j, m}/ (2 \pi i)$, we have
\begin{equation}\label{lem-E-gamma}
\gamma_{1, m} = - \frac{2 k_m + l_m}{3}, \quad \gamma_{2, m} = \frac{k_m -l_m}{3}, \quad \gamma_{3, m} = \frac{k_m+ 2 l_m}{3}. 
\end{equation}
It follows that 
\begin{equation*}
E_m =     \sum_{j=1}^3 \frac{\gamma_{j+1, m} - \gamma_{j, m}}{\gamma_{j+2, m}} 
=  \frac{3k_m}{k_m + 2 l_m } - \frac{3 l_m }{2 k_m + l_m} - \frac{3 (k_m+l_m)}{k_m -l_m}. 
\end{equation*}
Since 
\begin{multline*}
k_m(2k_m + l_m) (k_m - l_m) - l_m(k_m + 2 l_m) (k_m - l_m) - (k_m + l_m) (k_m + 2 l_m) (2 k_m + l_m) \\[6pt]
= 2 (k_m^2 - l_m^2) (k_m - l_m) - (k_m + l_m) (k_m + 2 l_m) (2 k_m + l_m) \\[6pt]
= (k_m + l_m) \Big( 2 k_m^2 - 4 k_m l_m + 2 k_m^2 - 2 k_m^2 - 2 l_m^2 - 5 k_m l_m \Big) = - 9 k_m l_m (k_m + l_m), 
\end{multline*}
we derive that 
\begin{equation*}
E_m
=  - \frac{2 7 k_m l_m (k_m+l_m)}{(k_m+2l_m) (2 k_m + l_m) (k_m-l_m)} \neq 0. 
\end{equation*}
The proof is complete. 
\end{proof}

We next show in \Cref{lem-2/3,lem-2pm} below that for  $L \in \cN$ and for  $1 \le m \le n_L$ with $p_m \neq 0$, it holds 
$$
2 p_m \neq 2 / (3 \sqrt{3}) \quad \mbox{ and } \quad 2 p_m \not \in \cP_L. 
$$
As a consequence $\varphi_{m, m}$ is constructed via 1) and 4) in \Cref{pro1}. We begin with 

\begin{lemma}\label{lem-2/3} Let $L \in \cN$ and $1 \le m \le n_L$.  Then 
$$
2 p_m \neq 2 / (3 \sqrt{3}). 
$$
\end{lemma}

\begin{proof} We first claim that there is no  $k, l \in \N_*$ with $k \ge l $ such that  
\begin{equation}\label{lem-2/3-p1}
(2k + l) (2 l + k) (k-l) = (k^2 + l^2 + kl)^{3/2}. 
\end{equation}
We prove this by contradiction. Assume that there exists such a pair $(k, l)$. Set
$$
H = \Big\{(k, l) \in  \N_*\times \N_*,  \; k \ge l, \mbox{ and \eqref{lem-2/3-p1} holds}\Big\}. 
$$
Set 
$$
h = \min \Big\{ k + l; (k, l) \in H \Big\} > 0.  
$$
Fix  $(k, l) \in H$ such that $k + l = h$. Since 
$$
(2k + l) (2 l + k) (k-l) \mbox{ is even}, 
$$
it follows from \eqref{lem-2/3-p1} that $k^2 + l^2 + kl $ is even. Hence both $k$ and $l$ are even. We write $k = 2 k_1$ and $l = 2 l_1$ for some $k_1, l_1 \in \N_*$. It is clear that 
$$
k_1 \ge l_1, 
$$
and
\begin{equation*}
(2k_1 + l_1) (2 l_1 + k_1) (k_1-l_1) = (k_1^2 + l_1^2 + k_1l_1)^{3/2}. 
\end{equation*}
This implies 
$$
(k_1, l_1) \in H. 
$$
We have
$$
k_1 + l_1 =  (k+ l)/ 2 = h/2 \quad \mbox{ and } \quad  h >  0.
$$
This contradicts the definition of $h$. The claim is proved. 

\medskip 
We are ready to derive the conclusion of \Cref{lem-2/3}. Since $2 p_m = 2 / (3 \sqrt{3})$ for some $1 \le m \le n_L$ and for some $L \in \cN$ if and only if, by the definition of $p_m$ in \eqref{def-pm},  
$$
(2k_m + l_m)(k_m - l_m)(2 l_m + k_m) = (k_m^2 + l_m^2 + k_m l_m)^{3/2}, 
$$
the conclusion follows from the claim. 
\end{proof}

We next prove 
\begin{lemma}\label{lem-2pm} There is no quadruple $(k_1, l_1, k_2, l_2) \in \N_*^4$ satisfying  the  system 
\begin{equation}\label{sys-kl}
\left\{\begin{array}{c}
k_1 > l_1,  \quad k_2 > l_2, \\[6pt]
k_1^2 + k_1 l_1 + l_1^2 = k_2^2 + k_2 l_2 + l_2^2,  \\[6pt]
(2k_2 + l_2) (2 l_2 + k_2) (k_2 - l_2) =  2 (2k_1 + l_1) (2 l_1 + k_1) (k_1 - l_1).  
\end{array}\right. 
\end{equation}
Consequently, for $L \in \cN$ and $1 \le m \le n_L$, we have 
\begin{equation}\label{lem-2pm-cl2}
2 p_m \not \in \cP_L \mbox{ if } p_m \neq 0. 
\end{equation}
\end{lemma}

\begin{proof} We prove the non-existence  by contradiction. Assume that there exists a quadruple $(k_1, l_1, k_2, l_2) \in \N_*^4$ satisfying \eqref{sys-kl}. 
Set 
\begin{equation}
G = \Big\{ (k_1, l_1, k_2, l_2) \in \N_*^4; \eqref{sys-kl} \mbox{ holds} \Big\}, 
\end{equation}
and let 
\begin{equation}
g = \min  \Big\{ k_1 + l_1 + k_2 + l_2;  (k_1, l_1, k_2, l_2) \in G  \Big\} > 0. 
\end{equation}
Fix $(k_1, l_1, k_2, l_2) \in G$ such that $k_1 + l_1 + k_2 + l_2 = g$. Set 
\begin{equation}\label{lem-2pm-A}
A := k_1^2 + k_1 l_1  + l_1^2 = k_2^2 + k_2 l_2  + l_2^2 \quad (\mbox{by the second line of \eqref{sys-kl}}). 
\end{equation}
Since, for $(k, l) \in \mR$,  
$$
(2k + l)(2l + k) = 2 (k^2 +  kl + l^2) + 3 kl \quad \mbox{ and } \quad (k-l)^2 = (k^2 + k l + l^2) - 3 k l, 
$$
it follows from the square of the last line of \eqref{sys-kl}, with 
\begin{equation}\label{def-x1x2}
x_1 = 3 k_1 l_1 \quad \mbox{ and } \quad x_2 = 3 k_2 l_2,
\end{equation}
that 
$$
(2 A + x_2)^2 (A - x_2) =  4 (2 A + x_1)^2 (A - x_1). 
$$
This implies 
\begin{equation}
(4 A^3 - 3A x_2^2 - x_2^3) = 4(4 A^3 - 3A x_1^2 - x_1^3), 
\end{equation}
or equivalently 
\begin{equation}\label{lem-2pm-p1}
12 A^3 = 3 A (4 x_1^2 - x_2^2) + 4 x_1^3 - x_2^3. 
\end{equation}
Using \eqref{def-x1x2}, we derive that $A^3 = 0 \mod 3$, which yields 
$$
A = 0 \mod 3. 
$$
Putting this information into \eqref{lem-2pm-p1} and using again \eqref{def-x1x2}, we obtain 
$$
x_1^3 - x_2^3  = 0 \mod 3^4. 
$$
We deduce from \eqref{def-x1x2} that 
\begin{equation}\label{lem-2pm-p2}
(k_1 l_1)^3 - (k_2 l_2)^3 = 0 \mod 3. 
\end{equation}
By writing $k_1 l_1$ under the form $k_2 l_2 + 3 q + r$ with $q \in \mZ$ and $r \in \N$ with $0 \le r \le 2$, we have 
\begin{equation}\label{lem-2pm-p3}
(k_1 l_1)^3 - (k_2 l_2)^3  = 3 k_2^2 l_2^2 (3 q + r) + 3 k_2 l_2 (3 q + r)^2 +  (3 q + r)^3. 
\end{equation}
Combining \eqref{lem-2pm-p2} and \eqref{lem-2pm-p3}  yields that $r = 0$. Putting this information into \eqref{lem-2pm-p1}, we obtain 
$$
A^3 = 0 \mod 3^4. 
$$
This implies 
$$
A = 0 \mod 9. 
$$
We deduce from \eqref{lem-2pm-A} that 
$$
k_1 = 0 \mod 3, \quad l_1 = 0 \mod 3, \quad k_2 = 0 \mod 3, \quad l_2 = 0 \mod 3. 
$$
Let $\hat k_1, \hat l_1, \hat k_2, \hat l_2 \in \N_*$ be such that 
$$
k_1 = 3 \hat k_1, \quad l_1 = 3 \hat l_1, \quad k_2 = 3 \hat k_2, \quad l_2 = 3 \hat l_2. 
$$
One can easily check that $(\hat k_1, \hat l_1, \hat k_2, \hat l_2) \in G$ and 
$$
\hat k_1 +  \hat l_1 +  \hat k_2 +  \hat l_2 = g/3 < g. 
$$
We obtain a contradiction. The non-existence associated with  \eqref{sys-kl} is proved. 

\medskip 
It is clear that \eqref{lem-2pm-cl2} is just a consequence of the non-existence by the definition of $L$  and $p_m$ as a function of $k_m$ and $l_m$ in  \eqref{def-L} and \eqref{def-pm}.   The proof is complete. 
\end{proof}

We are ready to state and prove the main result of this section:  

\begin{proposition}\label{pro-chi} Let $L \in \cN$ and $1 \le m \le n_L$. Then 
\begin{equation}\label{pro-chi-cl1}
\varphi_{m, m}' (0)  = 4 \pi L  =  - \phi_{m, m}'(0) \mbox{ if } p_m = 0, 
\end{equation}
and, if $p_m \neq 0$ and $s_m \neq 0$ then  
\begin{equation}\label{pro-chi-cl2}
\varphi_{m, m}' (0)   \neq 0. 
\end{equation}
\end{proposition}

\begin{proof} Assertion \eqref{pro-chi-cl1}  follows immediately from 4) of \Cref{pro1,pro1-Co}. We next consider the case $p_m \neq 0$. By \Cref{lem-2/3,lem-2pm}, 
we have 
$$
\varphi_{m, m}'(0) = 0 
$$
only if, with $\alpha = e^{2 \eta_{2, m} L}$ and $\lambda_j = \lambda_j (2 p_m)$,  
\begin{equation}
\left\{\begin{array}{c}
\sum_{j=1}^3 \lambda_j a_j = 0  \quad (= \varphi_{m, m}' (0) \mbox{ since $\chi_{m, m}'(0) = 0$}), \\[6pt]
\sum_{j=1}^3 \lambda_j e^{\lambda_j L} a_j = 0 \quad (= \varphi_{m, m}'(L)  \mbox{ since $\chi_{m, m}'(L) = 0$} ), \\[6pt]
\sum_{j=1}^3 (e^{\lambda_j L} - \alpha ) a_j = 0 \quad (= - \chi_{m, m}(L) + \alpha \chi_{m, m}(0) \mbox{ since $\chi_{m, m}(L) = \alpha \chi_{m, m } (0)$} ). 
\end{array}\right. 
\end{equation}
Since $E_m \neq 0$ by \Cref{lemE}, one has a non-trivial solution $(a_1, a_2, a_3)$ of this system. This implies 
\begin{equation}\label{pro-chi-p1}
\det K_1 = 0 \quad \mbox{ where } K_1 : = \quad \left(\begin{array}{ccc}
\lambda_1 & \lambda_2 & \lambda_3  \\[6pt]
\lambda_1 e^{\lambda_1 L} & \lambda_2 e^{\lambda_2 L}  & 
\lambda_3 e^{\lambda_3 L} \\[6pt]
e^{\lambda_1 L}  - \alpha & e^{\lambda_2 L}  - \alpha  & 
e^{\lambda_3 L}  - \alpha  
\end{array}\right). 
\end{equation}
Set 
$$
\hlambda_j = \lambda_j L. 
$$
Condition \eqref{pro-chi-p1} is equivalent to 
\begin{equation}\label{pro-chi-p2}
\det K_2 = 0 \quad \mbox{ where } \quad K_2 : = \left(\begin{array}{ccc}
\hlambda_1 & \hlambda_2 & \hlambda_3  \\[6pt]
\hlambda_1 e^{\hlambda_1} & \hlambda_2 e^{\hlambda_2}  & 
\hlambda_3 e^{\hlambda_3} \\[6pt]
e^{\hlambda_1}  - \alpha & e^{\hlambda_2}  - \alpha  & 
e^{\hlambda_3}  - \alpha  
\end{array}\right). 
\end{equation}
A computation yields 
\begin{equation*}
\det K_2 =  \sum_{j=1}^3 \hlambda_j \Big( ( \hlambda_{j+1} -  \hlambda_{j+2}) e^{\hlambda_{j+1} + \hlambda_{j+2}} - \alpha (\hlambda_{j+1}e^{\hlambda_{j+1}}  - \hlambda_{j+2}e^{\hlambda_{j+2}} ) \Big),  
\end{equation*}
which implies 
\begin{equation}\label{det-K2}
\det K_2 =   \sum_{j=1}^3 \hlambda_j  ( \hlambda_{j+1} -  \hlambda_{j+2}) \Big( e^{- \hlambda_{j}}  + \alpha e^{\hlambda_{j}} \Big). 
\end{equation}
Here we used the fact $\sum_{j=1}^3 \hlambda_j = L \sum_{j=1}^3 \lambda_j =0 $.  From the definition of $\lambda_j = \lambda_{j} (2 p_m)$ given in \Cref{def1}, we have 
\begin{equation*}
\left\{\begin{array}{c}
 \hlambda_1 + \hlambda_2 + \hlambda_3  = 0, \\[6pt]
\hlambda_1 \hlambda_2+ \hlambda_1 \hlambda_3 + \hlambda_2 \hlambda_3  = L^2, \\[6pt]
\hlambda_1 \hlambda_2 \hlambda_3 = 2 i p_m L^3. 
\end{array}\right. 
\end{equation*}
Define $\sigma_{j, m}$ by  
$$
\hlambda_{j} = \frac{2 \pi i  \sigma_{j, m}}{3}. 
$$
We then have 
\begin{equation*}
\left\{\begin{array}{c}
 \sigma_{1, m} +  \sigma_{2, m} +  \sigma_{3, m}  = 0, \\[6pt]
 \sigma_{1, m} \sigma_{2, m} + \sigma_{1, m} \sigma_{3, m} + \sigma_{2, m} \sigma_{3, m}  = - 3 (k_m^2 + l_m^2 + k_m l_m), \\[6pt]
\sigma_{1, m} \sigma_{2, m} \sigma_{3, m} = - 2 (2k_m + l_m) (2 l_m + k_m) (k_m - l_m),
\end{array}\right. 
\end{equation*}
where in the last identity, we used the fact
$$
p_m L^3 = \frac{1}{27} (2 \pi)^3 (2k_m + l_m) (2 l_m + k_m) (k_m - l_m). 
$$
It is clear that $\det K_2 = 0$ if and only if \eqref{main-assumption} holds. The proof is complete. 
\end{proof}

\section{Useful properties related to quasi-periodic functions}

In this section, we derive some properties for  $W_x(\cdot, 0)$ given in the introduction using the quasi-periodic-function theory. 
The main result of this section is \Cref{pro-quasi}. 
We begin with its weaker  version. 

\begin{lemma}\label{lem-quasi} Let  $\ell \in \N_*$, $a_j \in \mC$, $q_j  \ge 0$ for $1 \le j \le \ell$, and 
$M_{j_1, j_2}, N_{j_1, j_2} \in \mC$  with $1 \le j_1, j_2 \le \ell$. 
Assume that 
\begin{equation}\label{lem-quasi-A1}
\left\{\begin{array}{c}
q_{j_1} \neq q_{j_2} \mbox{ for } 1 \le j_1 \neq j_2 \le \ell, \\[6pt]
M_{j, j} \neq 0 \mbox{ for } 1 \le j  \le \ell, \\[6pt]
(M_{j, j} \mbox{ is real  and } N_{j, j} \neq 0)  \mbox{ if } q_j = 0, \\[6pt]
a_{j} \in i \mR \mbox{ if } q_j =0, 
\end{array}\right. 
\end{equation}
and
\begin{equation}\label{lem-quasi-A2}
\sum_{j=1}^\ell |a_{j}|^2  >  0. 
\end{equation}
Set, for $t \in \mR$, 
\begin{multline}\label{lem-quasi-g}
g(t) \\[6pt]:= \sum_{j_1 = 1}^\ell \sum_{j_2 = 1}^\ell \Big( a_{j_1} a_{j_2} M_{j_1, j_2} e^{- i (q_{j_1} + q_{j_2}) t} + \bar a_{j_1} \bar a_{j_2} \bar M_{j_1, j_2} e^{ i (q_{j_1} + q_{j_2}) t} + 2 a_{j_1} \bar a_{j_2} N_{j_1, j_2} e^{-i (q_{j_1} - q_{j_2})}  \Big). 
\end{multline}
There exists $t \in \mR_+$ such that 
\begin{equation}\label{lem-quasi-cl}
g(t)  \neq 0. 
\end{equation}
\end{lemma}

\begin{proof} We prove \eqref{lem-quasi-cl} by recurrence in $\ell$. It is clear that the conclusion holds for $\ell=1$. Indeed, if $q_1 \neq 0$ then since $e^{2 q_1 t}$, $0$, and $e^{-2 q_1 t}$ are independent, the conclusion follows. Otherwise, $q_1 = 0$. 
Since $M_{1, 1}$ is real and $a_{1} \in i \mR$, we have 
$$
g(t) =  2 |a_1|^2 N_{1, 1}. 
$$ 
The conclusion in the case $\ell = 1$ follows since $N_{1, 1} \neq 0$.

Assume that the conclusion holds for $\ell \ge 1$, we prove that the conclusion holds for $\ell +1$. 
Without loss of generality, one might assume that 
\begin{equation}\label{lem-quasi-q}
0 \le q_1 < q_2 < \dots < q_{\ell}  < q_{\ell +1}. 
\end{equation}
We will prove \eqref{lem-quasi-cl} for $\ell+1$ by contradiction. Assume that there exist $a_j$ and  $q_j \ge 0$ with $1 \le j \le \ell +1$, $M_{j_1, j_2}, \, N_{j_1, j_2} \in \mC$ with $1 \le j_1, j_2 \le \ell+1$ such that \eqref{lem-quasi-A1}, \eqref{lem-quasi-A2}, and \eqref{lem-quasi-q} hold, and,  for all $t  \in \mR_+$,  
\begin{equation}\label{lem-quasi-p1}
\sum_{j_1 = 1}^{\ell+1} \sum_{j_2 = 1}^{\ell + 1} \Big( a_{j_1} a_{j_2} M_{j_1, j_2} e^{- i (q_{j_1} + q_{j_2}) t} + \bar a_{j_1} \bar a_{j_2} \bar M_{j_1, j_2} e^{ i (q_{j_1} + q_{j_2}) t} + 2 a_{j_1} \bar a_{j_2} N_{j_1, j_2} e^{-i (q_{j_1} - q_{j_2})}  \Big)  = 0. 
\end{equation}
Since the function $e^{- 2 i q_{\ell+1} t}$ defined in $\mR_+$ does not belong to the space 
\begin{multline*}
\mbox{span} \left(  \Big \{e^{-it (q_{j_1} + q_{j_2})}; 1 \le j_1 \le \ell +1; 1 \le j_2 \le \ell \Big\},  \right. \\[6pt]
\left. \Big \{e^{it (q_{j_1} + q_{j_2})}; 1 \le j_1 \le \ell +1; 1 \le j_2 \le \ell + 1 \Big\},  \right. \\[6pt]
\left. \Big \{e^{-it (q_{j_1} - q_{j_2})}; 1 \le j_1 \le \ell +1; 1 \le j_2 \le \ell +1 \Big\}
\right), 
\end{multline*}
for $t \in \mR_+$ by  \eqref{lem-quasi-q},  we have 
$$
a_{\ell+1}^2 M_{\ell +1, \ell +1} = 0. 
$$
This yields, since $M_{\ell+1, \ell+1} \neq 0$, 
$$
a_{\ell+1} = 0. 
$$
It follows from \eqref{lem-quasi-p1} that 
\begin{equation}
\sum_{j_1 = 1}^{\ell} \sum_{j_2 = 1}^{\ell} \Big( a_{j_1} a_{j_2} M_{j_1, j_2} e^{- i (q_{j_1} + q_{j_2}) t} + \bar a_{j_1} \bar a_{j_2} \bar M_{j_1, j_2} e^{ i (q_{j_1} + q_{j_2}) t} + 2 a_{j_1} \bar a_{j_2} N_{j_1, j_2} e^{-i (q_{j_1} - q_{j_2})}  \Big)  = 0.  
\end{equation}
We now can use the assumption on the recurrence to obtain a contradiction. The proof of \eqref{lem-quasi-cl} is complete.  
\end{proof}

Using  \Cref{lem-quasi} and the theory of quasi-periodic functions, see e.g. \cite{Bohr47}, we can derive the following useful result for the proof of \Cref{thm1}. 

\begin{proposition}\label{pro-quasi} 

Let  $\ell \in \N_*$,  $a_j \in \mC$,  $q_j  \ge 0$ for $1 \le j \le \ell$, and 
$M_{j_1, j_2}, N_{j_1, j_2} \in \mC$  with $1 \le j_1, j_2 \le \ell$. 
Assume that \eqref{lem-quasi-A1} holds and denote $g$ by \eqref{lem-quasi-g}.  For all $0 < \gamma_1 < \gamma_2 $ there exist $\gamma_0>0$ and $\tau_0 > 0$ depending only on $\gamma_1$, $\gamma_2$,   $\ell$, $q_j$, $M_{j_1, j_2}$,  and $N_{j_1, j_2}$  such that if 
\begin{equation}
\gamma_1 \le  \sum_{j=1}^\ell |a_j|^2 \le  \gamma_2,
\end{equation}
then 
\begin{equation}\label{pro-cl}
\| g \|_{L^2(\tau, 2 \tau)} \ge \gamma_0 \mbox{ for all } \tau \ge \tau_0. 
\end{equation}
\end{proposition}

\begin{proof} Instead of \eqref{pro-cl}, it suffices to prove 
\begin{equation}\label{pro-cl1}
\| g \|_{L^\infty(\tau, 2 \tau)} \ge \gamma_0 \mbox{ for } \tau \ge \tau_0
\end{equation}
by contradiction since $|g'(t)| \le C$ in $\mR$.  
Assume that for all  $n \in \N_*$ there exist  $(a_{j, n})_{j=1}^\ell \subset \mC$ and $(t_n) \subset \mR$ such that $\gamma_1 \le  \sum_{j=1}^\ell |a_{j, n}|^2 \le  \gamma_2$, $t_n \ge n$,  and 
\begin{equation}\label{pro-p1}
\|g_n\|_{L^\infty(t_n, 2t_n)} \le 1/n, 
\end{equation}
where $g_n$ is defined in \eqref{lem-quasi-g} where $a_{j_1}$ and $a_{j_2}$ are replaced by  $a_{j_1, n}$ and $a_{j_2, n}$. Without loss of generality, one might assume that 
$$
\lim_{n \to + \infty} a_{j, n} = a_j \in \mC
$$
and $\gamma_1 \le  \sum_{j=1}^N |a_j|^2 \le  \gamma_2$. Consider $g$ defined by \eqref{lem-quasi-g} with these $a_j$. We have
\begin{equation}\label{pro-p2}
\lim_{n \to + \infty} \| g_n - g\|_{L^\infty(\mR)} = 0. 
\end{equation}

Since $g$ is an almost-periodic function with respect to $t$ (see e.g. \cite[Corollary on page 38]{Bohr47}), it follows from the definition of almost-periodic functions, see e.g. \cite[Section 44 on pages 32 and 33]{Bohr47}, that for every $\eps > 0$,  there exists $\cL_\eps > 0$
such that every interval $(\alpha, \alpha + \cL_\eps)$ containing a number $\tau(\eps, \alpha)$ for which it holds 
\begin{equation}\label{pro-ap}
|g(t + \tau(\eps, \alpha) ) - g(t)| \le \eps \mbox{ for all } t \in \mR. 
\end{equation}

The proof is now divided into two cases. 

\medskip 
\noindent{\it Case 1:}  $\liminf_{\eps \to 0 } \cL_\eps < + \infty$. Denote $\cL_0= \liminf_{\eps \to 0} \cL_\eps$. 
We claim that  $g$ is $T$-periodic for some period $T \le \cL_0 + 1$. Indeed, by \eqref{pro-ap} applied with $\alpha = 1/2$,  there exists a sequence $(\tau_n) \subset (1/2, \cL_0 +1)$ such that, for large $n$,  
$$
|g(t+ \tau_n) - g(t)| \le 1/n \mbox{ for all } t \in \mR. 
$$ 
By choosing $T = \liminf_{n \to + \infty} \tau_n$, we have 
$$
g(t + T) = g(t) \mbox{ for all } t \in \mR. 
$$
The claim is proved. 

Since $g$ is $T$-periodic, we have
$$
\|g \|_{L^\infty(t_n, t_n+  T + 1)} = \|g \|_{L^\infty(0, T + 1)}  \mbox{ for } n \in \N_*,
$$
and since $g$ is analytic and $g \neq 0$ by \Cref{lem-quasi}, we obtain 
$$
\|g \|_{L^\infty(0, T + 1)} > 0. 
$$
This contradicts \eqref{pro-p1} and \eqref{pro-p2}. The proof of Case 1 is complete. 

\medskip 
\noindent{\it Case 2:} $\lim_{\eps \to 0}L_\eps = + \infty$.  Set 
\begin{equation}\label{pro-quasi-rho}
\rho = \|g \|_{L^\infty(0, 1)}. 
\end{equation}
It follows from \Cref{lem-quasi} that $g$ is not identically equal to 0. Since $g$ is analytic, we derive  that  
\begin{equation}\label{pro-quasi-rho>0}
\rho > 0.
\end{equation}

Let $n_0 \ge 2$ be such that 
$$
\|g_n - g \|_{L^\infty(\mR)} < \rho/4, \quad \| g_n\|_{L^\infty(t_n, 2 t_n)} < \rho / 4 \mbox{ for } n \ge n_0. 
$$
Such an $n_0$ exists by \eqref{pro-p1}, \eqref{pro-p2}, and \eqref{pro-quasi-rho>0}.  We  have, for $n \ge n_0$,  
\begin{equation}\label{pro-p3}
\|g \|_{L^\infty(t_n, 2 t_n)} \le \|g_n - g\|_{L^\infty(t_n, 2 t_n)} + \|g_n\|_{L^\infty(t_n, 2 t_n)} \le \rho/4 + \rho/4 = \rho/2.  
\end{equation}
Fix $0< \eps < \rho/4$  and fix $n \ge n_0$ such that $1 \le \cL_\eps \le t_n/2$. Such a number $n$ exists since $t_n \ge n$. 
It follows from the definition of $\tau( \eps, t_n)$ that 
\begin{equation}\label{pro-tau1}
\tau(\eps, t_n) \in (t_n, t_n + \cL_\eps) \subset (t_n, 3 t_n/2), 
\end{equation}
and
\begin{equation}\label{pro-tau2}
\big|g\big(t + \tau(\eps, t_n)\big) - g(t)\big| \le \eps \mbox{ for all } t \in \mR. 
\end{equation}
This yields 
\begin{equation}\label{pro-p4}
\|g \|_{L^\infty(t_n, 2 t_n)} \mathop{\ge}^{\eqref{pro-tau1}} 
\|g \|_{L^\infty(\tau(\eps, t_n),  \tau(\eps, t_n) + 1)} \mathop{\ge}^{\eqref{pro-tau2}} 
\|g \|_{L^\infty(0, 1)} - \eps \ge \rho - \rho/4  = 3 \rho/4. 
\end{equation}
Combining \eqref{pro-p3} and \eqref{pro-p4} yields a contradiction since $\rho > 0$ by \eqref{pro-quasi-rho>0}. The proof of Case 2 is complete. 
\end{proof}

\section{An upper bound for the decay rate - Proof of \Cref{thm1}} \label{sect-thm1}

This section containing two subsections is devoted to the proof of \Cref{thm1}. The main ingredient is given in the first section and the proof is presented in the second one. 

\subsection{A key lemma} In this section, we prove 

\begin{lemma}\label{lemK} Let $L \in \cN$. Assume that $\dim \M = 1$ or  \eqref{main-assumption} holds.  There exist $\eps_0 > 0$, $C>0$,  and $T_0 > 0$ depending only on $L$ such that for all (real) $u_0 \in L^2(0, L)$ with $\| u_0 \|_{L^2(0, L)} \le \eps_0$, the unique solution $u \in C \big([0, + \infty); L^2(0, L) \big) \cap L^2_{\loc} \big([0, + \infty); H^1(0, L) \big)$
of system \eqref{KdV-NL}  satisfies 
\begin{equation}\label{decayK}
\| u(T, \cdot) \|_{L^2(0, L)} \le  \| u_0 \|_{L^2(0, L)} \Big(1 -  C \| u_0 \|_{L^2(0, L)}^2 \Big) \mbox{ for } T \ge T_0. 
\end{equation}
\end{lemma}

\begin{proof} We first collect several known facts. Let $T_1 > 0$ be such that 
\begin{equation}\label{lemK-T1}
\| v_x(\cdot, 0) \|_{L^2(0, t)} \ge \frac{1}{2} \| v(0, \cdot) \|_{L^2(0, L)} \mbox{ for } t \ge T_1, 
\end{equation}
for all solutions $v \in C \big([0, + \infty); L^2(0, L) \big) \cap L^2_{\loc} \big([0, + \infty); H^1(0, L) \big)$  of the system 
\begin{equation}\left\{
\begin{array}{cl}
v_t (t, x) + v_x (t, x) + v_{xxx} (t, x)   = 0 &  \mbox{ in } (0, +\infty)  \times  (0, L), \\[6pt]
v(t, x=0) = v(t, x=L) = v_x(t , x= L)=  0 & \mbox{ in }  (0, +\infty), 
\end{array}\right.
\end{equation}
with $v(0, \cdot) \in L^2(0, L)$ satisfying the condition 
$$
v(0, \cdot) \perp \M 
$$
(the orthogonality is considered with respect to $L^2(0, L)$-scalar product).  The existence of such a constant $T_1$ follows from \cite{Rosier97}.

There exist two positive constants $\eps_0$ and  $c_1$ such that if $\| u_0 \|_{L^2(0, L)} \le \eps_0$, then 
\begin{equation}\label{lemK-c1}
\| u \|_{C\big([0, T_1]; L^2(0, L) \big)} + \| u \|_{L^2\big( (0, T_1); H^1(0, L) \big)} \le c_1 \| u_0\|_{L^2(0, L)}
\end{equation}
(see e.g., \cite[Proposition 14]{CC04}). 

There is  a positive constant $c_2$ such that if $\tu_0 \in L^2(0, L)$, $\tf \in L^1\big((0, T_1); L^2(0, L) \big)$,  and $\ty \in C \big([0, T_1); L^2(0, L) \big) \cap L^2 \big([0, T_1); H^1(0, L) \big) $ is the unique solution of the system 
\begin{equation}\left\{
\begin{array}{cl}
\tu_t (t, x) + \tu_x (t, x) + \tu_{xxx} (t, x)   = \tf &  \mbox{ in } (0, T_1) \times (0, L), \\[6pt]
\tu(t, x=0) = \tu(t, x=L) = \tu_x(t , x= L)=  0 & \mbox{ in }  (0, T_1), \\[6pt]
\tu(t= 0, \cdot) = \tu_0 & \mbox{ in } (0, L), 
\end{array}\right.
\end{equation}
then 
\begin{multline}\label{lemK-c2}
\| \tu_x(\cdot, 0) \|_{L^2(0, T_1)} + \| \tu \|_{C\big([0, T_1]; L^2(0, L) \big)} + \| \tu \|_{L^2\big( (0, T_1); H^1(0, L) \big)} \\[6pt]
 \le c_2 \Big( \| \tu_0\|_{L^2(0, L)} + \|\tf\|_{L^1\big((0, T_1); L^2(0, L) \big)} \Big). 
\end{multline}

There exists a positive constant $c_3$ depending only on $L$ such that, for {\it all} $T>0$, 
\begin{equation}\label{lemK-c3} 
\| \xi \xi_x \|_{L^1\big( (0, T); L^2(0, L) \big)} \le c_3 \| \xi \|_{L^2\big( (0, T); H^1(0, L) \big)}^2 
\end{equation}
(the constant $c_3$ is independent of $T$). 

\medskip 
We now decompose $u_0$ into two parts: 
\begin{equation}\label{thm-u0-decomp}
u_0 = u_{0, 1} + u_{0, 2} \mbox{ in }  (0, L), 
\end{equation}
where
\begin{equation*}
u_{0, 1} = \mbox{Projection}_{\M} u_0  
\end{equation*}
with respect to $L^2(0, L)$-scalar product.

\medskip 
The proof is now divided  into two cases, with $0< \eps = \| u_{0} \|_{L^2(0, L)} < \eps_0$ (the conclusion is clear if $\eps =0$), 
\begin{itemize}
\item Case 1:   $\| u_{0, 2} \|_{L^2(0, L)} \ge \beta \eps^2 = \beta \| u_0\|_{L^2(0, L)}^2$,  
\item Case 2: $\| u_{0, 2} \|_{L^2(0, L)} <  \beta \eps^2 = \beta \| u_0\|_{L^2(0, L)}^2$, 
\end{itemize}
where 
\begin{equation}\label{def-beta}
\beta = 4c_1^2 c_2 c_3.
\end{equation}

\medskip
\noindent{\it Case 1}: Assume that 
\begin{equation}\label{Case1}
\| u_{0, 2} \|_{L^2(0, L)} \ge \beta \eps^2 = \beta \| u_0\|_{L^2(0, L)}^2. 
\end{equation}
Let $\hu \in C \big([0, T_1); L^2(0, L) \big) \cap L^2 \big([0, T_1); H^1(0, L) \big)$ be the unique solution of 
\begin{equation}\label{lemK-hu}\left\{
\begin{array}{cl}
\hu_t (t, x) + \hu_x (t, x) + \hu_{xxx} (t, x) = 0 &  \mbox{ in } (0, T_1) \times (0, L), \\[6pt]
\hu(t, 0) = \hu(t, L) = \hu_x(t , L)=  0 & \mbox{ in } (0, T_1), \\[6pt]
\hu(0, \cdot)  = u_0 & \mbox{ in } (0, L). 
\end{array}\right.
\end{equation}
Then 
\begin{equation}\label{lemK-c1c2}
\|(\hu - u)_x(\cdot, 0) \|_{L^2(0, T_1)} \mathop{\le}^{\eqref{lemK-c2}} c_2 \| u u_x \|_{L^1 \big((0, T_1) ; L^2(0, L) \big)} \mathop{\le}^{\eqref{lemK-c1}, \eqref{lemK-c3}} c_1^2 c_2 c_3  \eps^2. 
\end{equation}

Let $\hu_j \in C \big([0, T_1); L^2(0, L) \big) \cap L^2\big([0, T_1); H^1(0, L) \big)$ with $j=1, \, 2$ be the unique solution of  
\begin{equation}\left\{
\begin{array}{cl}
\hu_{j, t} (t, x) + \hu_{j, x} (t, x) + \hu_{j, xxx} (t, x) = 0 &  \mbox{ for } t \in (0, T), \, x \in (0, L), \\[6pt]
\hu_j(t, 0) = \hu_j(t, L) = \hu_{j, x} (t , L)=  0 & \mbox{ for } t \in (0, T), \\[6pt]
\hu_j(0, \cdot)  = u_{0, j} &  \mbox{ in } (0, L). 
\end{array}\right.
\end{equation}
Then  
$$
\hu = \hu_1 + \hu_2 \mbox{ in } [0, T_1] \times [0, L]. 
$$
We have
\begin{equation}\label{Case1-p1}
\hu_{1, x} (\cdot, 0) =0 \mbox{ in } [0, T_1],    
\end{equation}
and,  by the choice of $T_1$ via \eqref{lemK-T1}, 
\begin{equation}\label{Case1-p2}
\| \hu_{2, x} (\cdot, 0)\|_{L^2(0, T_1)} \ge \frac{1}{2}  \| \hu_2(0, \cdot)\|_{L^2(0, L)} = \frac{1}{2}  \| u_{0,2}\|_{L^2(0, L)}. 
\end{equation}
It follows from \eqref{Case1} that 
\begin{equation}\label{Case1-p3}
\| \hu_{x} (\cdot, 0)\|_{L^2(0, T_1)} \ge \frac{1}{2} \beta \eps^2. 
\end{equation}

From \eqref{lemK-c1c2} and \eqref{Case1-p3}, we obtain
\begin{multline*}
\| u_{x} (\cdot, 0)\|_{L^2(0, T_1)}  \ge  \| \hu_{x} (\cdot, 0)\|_{L^2(0, T_1)} - \| (u - \hu)_{x} (\cdot, 0)\|_{L^2(0, T_1)}  \ge \left( \frac{1}{2} \beta - c_1^2 c_2 c_3 \right) \eps^2 \mathop{\ge}^{\eqref{def-beta}} c_1^2 c_2 c_3 \eps^2. 
\end{multline*}
In other words,  
\begin{equation}\label{Case1-*}
\| u_{x} (\cdot, 0)\|_{L^2(0, T_1)} 
 \ge  c_1^2 c_2 c_3 \| u_0\|_{L^2(0, L)}^2. 
\end{equation}

\medskip 
\noindent{\it Case 2}: Assume that 
\begin{equation}\label{Case2}
\| u_{0, 2} \|_{L^2(0, L)} < \beta \eps^2 =  \beta \| u_0 \|_{L^2(0, L)}^2. 
\end{equation}
Since 
$$
\| u_{0, 1} \|_{L^2(0, L)}^2 + \| u_{0, 2} \|_{L^2(0, L)}^2 = \| u_{0} \|_{L^2(0, L)}^2 = \eps^2, 
$$
by considering $\eps$ sufficiently small, one can assume that
$$
\| u_{0, 1} \|_{L^2(0, L)} \ge \eps/2. 
$$ 

Let $\alpha_m \in \mC$ ($1 \le m \le n_L$) be such that 
\begin{equation}\label{thm-def-y1}
\frac{1}{\eps} u_{0, 1} = \Re \left\{ \sum_{m=1}^{n_L}  \alpha_m  \Psi_m (0, x) \right\}. 
\end{equation}
Since $u_{0, 1} \in \M$, such a family of $(\alpha_m)_{m=1}^{n_L}$ exists. Since 
$$
1/2 \le \| \frac{1}{\eps} u_{0, 1}  \|_{L^2(0, L)} \le 1
$$
and $\Big( \Psi_m (0, \cdot) \Big)$ is orthogonal in $L^2(0, L)$ (with respect to the complex field), one can assume in addition that 
$$
0< \gamma_1 \le \sum_{m=1}^{n_L}  |\alpha_m|^2 \le \gamma_2, 
$$
for some constants $\gamma_1$, $\gamma_2$ depending only on $L$. Moreover, since $\Psi_m (0, x)  \in i \mR$ for $x \in [0, L]$ (by \eqref{psi-complex}) if $p_m = 0$ (see e.g. \eqref{psi-complex}), one can also assume that $a_m \in i \mR$ if $p_m =0$.

Let $\gamma_0 > 0$  and  $\tau_0> 0$ be  the constants given  in  \Cref{pro-quasi} with $\ell= n_L$, $\gamma_1$ and  $\gamma_2$ determined above, $q_m= p_m$  given by \eqref{def-pm}, 
\begin{equation}\label{Case2-M}
M_{m_1, m_2} =  \frac{1}{8} \varphi_{m_1, m_2}' (0) \quad \mbox{ and } \quad N_{m_1, m_2} = \frac{1}{8}  \phi_{m_1, m_2}' (0), 
\end{equation}
where $\varphi_{m_1, m_2}$ and $\phi_{m_1, m_2}$ are defined in \Cref{pro1} and \Cref{pro1-Co}, respectively; in the case the definition of $\varphi_{m_1, m_2}$ and $\phi_{m_1, m_2}$ in \Cref{pro1} and \Cref{pro1-Co} are not unique, we fix  a choice of  $\varphi_{m_1, m_2}$ and $\phi_{m_1, m_2}$.

By \Cref{pro-chi}, we have 
$$
M_{m,m} \neq 0, 
$$
and 
$$
(M_{m, m}  \mbox{ is real and $N_{m, m} \neq 0$) if } p_m = 0. 
$$

Then, by \Cref{pro-quasi}, for all 
$a_j \in \mC$ $(1 \le j \le N)$ satisfying 
$\gamma_1 \le  \sum_{j=1}^N |a_j|^2 \le \gamma_2$,  it holds 
\begin{equation}\label{Case2-pro-g}
\| g \|_{L^2(\tau, 2 \tau)} \ge \gamma_0 \mbox{ for all } \tau \ge \tau_0, 
\end{equation}
where 
\begin{multline}\label{Case2-def-g}
g(t) = \sum_{m_1 = 1}^{n_L} \sum_{m_2 = 1}^{n_L} \Big( a_{m_1} a_{m_2} M_{m_1, m_2} e^{- i (p_{m_1} + p_{m_2}) t} \\[6pt]  + \bar a_{m_1} \bar a_{m_2} \bar M_{m_1, m_2} e^{ i (p_{m_1} + p_{m_2}) t} 
+ 2 a_{m_1} \bar a_{m_2} N_{m_1, m_2} e^{-i (p_{m_1} - p_{m_2})}  \Big).
\end{multline}

Define
\begin{equation}\label{Case2-defA}
A  = \beta + 2  \sum_{m_1 = 1}^{n_L} \sum_{m_2 = 1}^{n_L}  \|\varphi_{m_1, m_2}\|_{L^2(0, L)} + 2 \sum_{m_1 = 1}^{n_L} \sum_{m_2 = 1}^{n_L}  \|\phi_{m_1, m_2}\|_{L^2(0, L)}, 
\end{equation}
and set 
\begin{equation}\label{Case2-def-c4}
c_4 = 1/ (2 A). 
\end{equation}

Let  $T_2 \ge 2 \tau_0$ be such that 
\begin{equation}\label{Case2-T2}
\| y_x(\cdot, 0) \|_{L^2(T_2/2, T_2)} \le c_4 \gamma_0  \| y(0, \cdot) \|_{L^2(0, L)}, 
\end{equation}
for all solutions  $y \in C \big([0, + \infty); L^2(0, L) \big) \cap L^2_{\loc} \big([0, + \infty); H^1(0, L) \big)$ of 
\begin{equation}\label{thm-sys1}\left\{
\begin{array}{cl}
y_t (t, x) + y_x (t, x) + y_{xxx} (t, x)   = 0 &  \mbox{ for } t \in (0, +\infty), \, x \in (0, L), \\[6pt]
y(t, 0) = y(t, L) = y_x(t , L)=  0 & \mbox{ for } t \in (0, +\infty), 
\end{array}\right.
\end{equation}
with $y(0, \cdot) \in L^2(0, L)$. Note that $T_2$ is independent of $y(0, \cdot)$. The existence of $T_2$ can be proved by decomposing $y(0, \cdot) = y_1(0, \cdot) + y_2(0, \cdot)$ with $y_1(0, \cdot) \in \M$, and noting that \eqref{Case2-T2} holds for the solution with initial data being $y_2(0, \cdot)$ 
since the solution is exponential decay,  and the contribution for $y_x(\cdot, 0)$ from  the solution with initial data is $y_1(0, \cdot)$ is 0.

Let $\tu_1, \; \tu_2 \in C \big([0, + \infty); L^2(0, L) \big) \cap L^2_{\loc} \big([0, + \infty); H^1(0, L) \big)$ be the unique solution of 
\begin{equation}\label{sys-NL}\left\{
\begin{array}{cl}
\tu_{1, t} (t, x) + \tu_{1, x} (t, x) + \tu_{1, xxx} (t, x)  = 0 &  \mbox{ for } t \in (0, +\infty), \, x \in (0, L), \\[6pt]
\tu_1(t, 0) = \tu_1(t, L) = \tu_{1, x}(t , L)=  0 & \mbox{ for } t \in (0, +\infty), \\[6pt]
\dsp \tu_1(0, \cdot ) =  \frac{1}{\eps} u_{0, 1} & \mbox{ in } [0, L], 
\end{array}\right.
\end{equation} 
and 
\begin{equation}\label{Case2-sys}\left\{
\begin{array}{cl}
\tu_{2, t} (t, x) + \tu_{2, x} (t, x) + \tu_{2, xxx} (t, x) + \tu_1 \tu_{1, x}  = 0 &  \mbox{ for } t \in (0, +\infty), \, x \in (0, L), \\[6pt]
\tu_2(t, 0) = \tu_2(t, L) = \tu_{2, x}(t , L)=  0 & \mbox{ for } t \in (0, +\infty), \\[6pt]
\dsp \tu_2(0, \cdot)  = \frac{1}{\eps^2} u_{0, 2} & \mbox{ in } [0, L]. 
\end{array}\right.
\end{equation}

Set 
$$
V(t, x) =  \sum_{m=1}^{n_L} \alpha_m \Psi_m(t, x) \quad \mbox{ and } \quad U(t, x) = \Re V(t, x). 
$$
We have 
\begin{equation}\label{sys-U}\left\{
\begin{array}{cl}
U(t, x) + U_x (t, x) + U_{xxx} (t, x)  = 0 &  \mbox{ for } t \in (0, +\infty), \, x \in (0, L), \\[6pt]
U(t, x=0) = U(t, x=L) = U_x(t , x= L)=  0 & \mbox{ for } t \in (0, +\infty), \\[6pt]
\dsp U(t=0, \cdot) = \frac{1}{\eps} u_{0, 1} & \mbox{ in } [0, L].  
\end{array}\right.
\end{equation}
This implies 
$$
\tu_1 = U \mbox{ in } (0, + \infty) \times (0, L). 
$$

Define 
\begin{equation}\label{def-V1}
V_1(t, x) = \sum_{m_1=1}^{n_L} \sum_{m_2=1}^{n_L} \alpha_{m_1} \alpha_{m_2} \varphi_{m_1, m_2}(x) e^{- i (p_{m_1}  + p_{m_2} ) t}, 
\end{equation}
and 
\begin{equation}\label{def-V2}
V_2(t, x) =  \sum_{m_1=1}^{n_L} \sum_{m_2=1}^{n_L} \alpha_{m_1} \bar \alpha_{m_2} \phi_{m_1, m_2}(x) e^{- i (p_{m_1}  - p_{m_2} ) t}.  
\end{equation}
Then, by the construction of $\varphi_{m_1, m_2}$,  
\begin{equation}\left\{
\begin{array}{cl}
V_{1, t} (t, x) + V_{1, x} (t, x) + V_{1, xxx} (t, x) + \big( V(t, x) V(t, x) \big)_x = 0 &  \mbox{ for } t \in (0, +\infty), \, x \in (0, L), \\[6pt]
V_1(t, 0) = V_1(t, L) = V_{1, x}(t , L)=  0 & \mbox{ for } t \in (0, +\infty), 
\end{array}\right.
\end{equation}
and, by the construction of $\phi_{m_1, m_2}$,  
\begin{equation}\left\{
\begin{array}{cl}
V_{2, t} (t, x) + V_{2, x} (t, x) + V_{2, xxx} (t, x) + \big( |V(t, x)|^2 \big)_x  = 0 &  \mbox{ for } t \in (0, +\infty), \, x \in (0, L), \\[6pt]
V_2(t, 0) = V_2(t, L) = V_{2, x}(t , L)=  0 & \mbox{ for } t \in (0, +\infty). 
\end{array}\right.
\end{equation}

Set 
\begin{equation}\label{def-W}
W = \frac{1}{8} \Big( V_1 + \bar V_1 + 2 V_2 \Big) \mbox{ in } (0, + \infty) \times (0, L). 
\end{equation}
It follows from  \eqref{Case2-def-g} that $W_x(t, 0) = g(t)$ in $\mR_+$ and hence, by \eqref{Case2-pro-g}, 
\begin{equation}\label{Case2-pro-W}
\| W_x(t, 0) \|_{L^2(\tau, 2 \tau)} \ge \gamma_0 \mbox{ for all } \tau \ge \tau_0.  
\end{equation}

Since 
$$
\big( V(t, x) V(t, x) \big)_x + \overline{\big( V(t, x) V(t, x) \big)_x} + 2 \big( |V(t, x)|^2 \big)_x = 8 U(t, x) U_x(t, x), 
$$
we derive from \eqref{def-W} that 
\begin{equation}\left\{
\begin{array}{cl}
W_{t} (t, x) + W_{x} (t, x) + W_{xxx} (t, x)  + U(t, x) U_{x}(t, x) = 0 &  \mbox{ for } t \in (0, +\infty), \, x \in (0, L), \\[6pt]
W(t, 0) = W(t, L) = W_{x}(t , L)=  0 & \mbox{ for } t \in (0, +\infty). \end{array}\right.
\end{equation}

Let $\tW \in C \big([0, + \infty); L^2(0, L) \big) \cap L^2_{\loc} \big([0, + \infty); H^1(0, L) \big)$ be the unique solution of 
\begin{equation}\left\{
\begin{array}{cl}
\tW_{t} (t, x) + \tW_{x} (t, x) + \tW_{xxx} (t, x)   = 0 &  \mbox{ for } t \in (0, +\infty), \, x \in (0, L), \\[6pt]
\tW(t, 0) = \tW(t, L) = \tW_{x}(t , L) =  0 & \mbox{ for } t \in (0, +\infty), \\[6pt]
\dsp \tW(0, \cdot) = \tu_{2} (0, \cdot) - W(0, \cdot). 
\end{array}\right.
\end{equation}
Then 
\begin{equation}\label{Case2-hu2WW}
\tu_2 = \tW + W \mbox{ in } (0, + \infty) \times (0, L). 
\end{equation}

We have 
\begin{equation}\label{Case2-tt1}
\| \tu_{2, x} (\cdot, 0) \|_{L^2(T_2/2, T_2)} \mathop{\ge}^{\eqref{Case2-hu2WW}} \| W_{x} (\cdot, 0) \|_{L^2(T_2/2, T_2)} - \| \tW_{x} (\cdot, 0) \|_{L^2(T_2/2, T_2)}, 
\end{equation}
\begin{equation}\label{Case2-tt2}
\| \tW_{x} (\cdot, 0) \|_{L^2(T_2/2, T_2)} \mathop{\le}^{\eqref{Case2-T2}} c_4 \gamma_0 \|\tW(0, \cdot) \|_{L^2(0, L)}, 
\end{equation}
and, since $T_2 \ge \tau_0$,  
\begin{equation}\label{Case2-tt3}
\| W_{x} (\cdot, 0) \|_{L^2(T_2/2, T_2)} \mathop{\ge}^{\eqref{Case2-pro-W}} \gamma_0. 
\end{equation}

Since, by \eqref{def-V1}, \eqref{def-V2}, and \eqref{def-W} 
\begin{multline*}
8W(0, x) = \sum_{m_1=1}^{n_L} \sum_{m_2=1}^{n_L} \alpha_{m_1} \alpha_{m_2} \varphi_{m_1, m_2}(x)  \\[6pt]
+ \sum_{m_1=1}^{n_L} \sum_{m_2=1}^{n_L} \bar \alpha_{m_1} \bar \alpha_{m_2}  \bar \varphi_{m_1, m_2}(x)
+  2 \sum_{m_1=1}^{n_L} \sum_{m_2=1}^{n_L} \alpha_{m_1} \bar \alpha_{m_2} \phi_{m_1, m_2}(x), 
\end{multline*}
it follows that 
\begin{equation}\label{Case2-W0}
\|W(0, \cdot) \|_{L^2(0, L)} \le 2  \sum_{m_1 = 1}^{n_L} \sum_{m_2 = 1}^{n_L}  \|\varphi_{m_1, m_2}\|_{L^2(0, L)} + 2 \sum_{m_1 = 1}^{n_L} \sum_{m_2 = 1}^{n_L}  \|\phi_{m_1, m_2}\|_{L^2(0, L)}. 
\end{equation}
By the definition of $A$ in  \eqref{Case2-defA}, we obtain   from \eqref{Case2} and \eqref{Case2-W0} that 
\begin{equation}\label{Case2-A2}
A  \ge 
\| \tu_2 (0, \cdot) \|_{L^2(0, L)} + \|W(0, \cdot) \|_{L^2(0, L)} \mathop{\ge}^{\eqref{Case2-hu2WW}} \| \tW(0, \cdot) \|_{L^2(0, L)}. 
\end{equation}

Combining \eqref{Case2-tt1},  \eqref{Case2-tt2}, \eqref{Case2-tt3}, and \eqref{Case2-A2} yields
\begin{equation*}
\|\tu_{2, x} (\cdot, 0) \|_{L^2(T_2/2, T_2)}  \ge \gamma_0  - c_4 \gamma_0 A. 
\end{equation*}
Since $c_4 = 1 / (2 A)$ by \eqref{Case2-def-c4},  we  obtain 
\begin{equation}\label{lemK-p1}
\|\tu_{2, x} (\cdot, 0) \|_{L^2(T_2/2, T_2)}  \ge \gamma_0/2. 
\end{equation}
  
Set 
\begin{equation*}
u_d = \eps \tu_1 + \eps^2 \tu_2 - u \mbox{ in } (0, + \infty) \times (0, L), 
\end{equation*}
and 
\begin{equation}
f_d = u u_x - \eps^2 \tu_1 \tu_{1, x} \mbox{ in } (0, + \infty) \times (0, L). 
\end{equation}
We have, by \eqref{sys-NL} and \eqref{Case2-sys}, 
\begin{equation}\left\{
\begin{array}{cl}
u_{d, t} (t, x) + u_{d, x} (t, x) + u_{d, xxx} (t, x)  = f_d (t, x) &  \mbox{ for } t \in (0, +\infty), \, x \in (0, L), \\[6pt]
u_d(t, x=0) = u_d(t, x=L) = u_d(t , x= L)=  0 & \mbox{ for } t \in (0, +\infty), \\[6pt]
u_d(t = 0, \cdot)  = 0 & \mbox{ in } (0, L). 
\end{array}\right.
\end{equation}
It is clear that 
\begin{equation*}
\| f_d\|_{L^1\big((0, T_2); L^2(0, L) \big)} \le C \eps^2, 
\end{equation*}
where $C$ is a positive constant depending only on  $T_2$ and $L$. It follows that 
\begin{equation*}
\| u_d \|_{C\big([0, T_2]; L^2(0, L) \big)} + \| u_d \|_{L^2\big( (0, T_2); H^1(0, L) \big)} \le C \eps^2. 
\end{equation*}
This in turn implies that 
\begin{equation*}
\| f_d\|_{L^1\big((0, T_2); L^2(0, L) \big)} \le C \eps^3. 
\end{equation*}
and therefore, 
\begin{equation}\label{lemK-p2}
\| u_{d, x} (\cdot, 0)  \|_{L^2( 0, T_2) } \le C \eps^3. 
\end{equation}

Combining \eqref{lemK-p1} and  \eqref{lemK-p2}, and noting that $\tu_{1, x} (t, 0) = 0$ yield 
\begin{equation}\label{est-case2}
\|u_x (\cdot, 0)  \|_{L^2(T_2/2, T_2) } \ge C \eps^2.  
\end{equation}
The analysis of Step 2 is complete. 

\medskip 

The conclusion now follows from Case 1 where one obtains \eqref{Case1-*} and Case 2 where one obtains  \eqref{est-case2} by choosing $T_0 = \max\{T_1, T_2\}$ and  using \eqref{key-identity}. 
The proof is complete. 
\end{proof}


We are ready to give 

\subsection{Proof of \Cref{thm1}} By \Cref{lemK}, we have 
$$
\|u(T_2, \cdot)  \|_{L^2(0, L) } \le \| u(0, \cdot) \|_{L^2(0, L)} \Big(1 - C \| u(0, \cdot) \|_{L^2(0, L)}^2 \Big).  
$$
This yields, with $\| u(0, \cdot) \|_{L^2(0, L)} = \eps>0$ and $p$ being the largest integer less than $1/ (2 C \eps^2)$, 
$$
\|u(p T_2, \cdot)  \|_{L^2(0, L) } \le \frac{1}{2} \| u(0, \cdot) \|_{L^2(0, L)}.  
$$
Here we also used \eqref{key-identity}. Using \eqref{key-identity} again, it follows that, for $T \ge C/ \| u(0, \cdot) \|_{L^2(0, L)}^2$,
$$
\|u(T, \cdot)  \|_{L^2(0, L) } \le \frac{1}{2} \| u(0, \cdot) \|_{L^2(0, L)}. 
$$
This implies, by recurrence, that  
$$
\| u(T, \cdot) \|_{L^2(0, L)} \le  2^{-n} \| u_0 \|_{L^2(0, L)}\mbox{ for } T \ge  C  \sum_{p=0}^{n-1} 2^{2p} /\| u(0, \cdot) \|_{L^2(0, L)}^2 
$$
since $\|u(t, \cdot)  \|_{L^2(0, L)}$ is a non-increasing function with respect to $t$. In particular, we obtain, since $\|u(t, \cdot)  \|_{L^2(0, L)}$ is a non-increasing function with respect to $t$, 
\begin{equation}
 \|u(t, \cdot)  \|_{L^2(0, L)} \le C/ t^{1/2}.  
\end{equation}
The proof is complete. \qed

\section{A lower bound for the decay rate - Proof of \Cref{pro-opt}} \label{sect-opt}

Fix $1 \le m \le n_L$ and  $\alpha_m \in \mC$ with $|\alpha_m| = 1$ such that 
$$
\Re ( \alpha_m \varphi_{m, m}(x)) \mbox{ is not identically equal to 0 in } [0, L]. 
$$
Let $\tu_1 \in C \big([0, + \infty); L^2(0, L) \big) \cap L^2_{\loc} \big([0, + \infty); H^1(0, L) \big)$ be the unique solution of 
\begin{equation}\label{sys-NL-L}\left\{
\begin{array}{cl}
\tu_{1, t} (t, x) + \tu_{1, x} (t, x) + \tu_{1, xxx} (t, x)  = 0 &  \mbox{ for } t \in (0, +\infty), \, x \in (0, L), \\[6pt]
\tu_1(t, x=0) = \tu_1(t, x=L) = \tu_{1, x}(t , x= L)=  0 & \mbox{ for } t \in (0, +\infty), \\[6pt]
\dsp \tu_1(0, \cdot ) = \Re ( \alpha_m \varphi_{m, m}).  
\end{array}\right.
\end{equation}

Set 
\begin{equation}\label{def-V1-L}
V_1(t, x) = \alpha_{m}^2  \varphi_{m, m}(x) e^{- 2 i p_{m} t}, 
\end{equation}
\begin{equation}\label{def-V2-L}
V_2(t, x) =  |\alpha_{m}|^2  \phi_{m, m}(x), 
\end{equation}
and denote  
\begin{equation}\label{def-W-L}
\tu_2 = \frac{1}{8} \Big( V_1 + \bar V_1 + 2 V_2 \Big) \mbox{ in } (0, + \infty) \times (0, L). 
\end{equation}
Since $\phi_{m, m}$ is real by 3) of \Cref{pro1-Co}, it follows that $V_2$ is real and hence so is $\tu_2$. 

As in the proof of \Cref{lemK}, we have
$$
\tu_1(t, x) = \Re \Big(\alpha_m \varphi_{m, m} (x)e^{-  i p_m t} \Big), 
$$
and 
\begin{equation*}\left\{
\begin{array}{cl}
\tu_{2, t} (t, x) + \tu_{2, x} (t, x) + \tu_{2, xxx} (t, x)  + \tu_1(t, x) \tu_{1, x}(t, x) = 0 &  \mbox{ for } t \in (0, +\infty), \, x \in (0, L), \\[6pt]
\tu_2(t, x=0) = \tu_2(t, x=L) = \tu_{2, x}(t , x= L)=  0 & \mbox{ for } t \in (0, +\infty). \end{array}\right.
\end{equation*}

Let $u \in C \big([0, + \infty); L^2(0, L) \big) \cap L^2_{\loc} \big([0, + \infty); H^1(0, L) \big)$ be a (real) solution of \eqref{KdV-NL} with 
$$
\| u(0, \cdot) \|_{L^2(0, L)} \le \Gamma \eps, 
$$
where 
$$
\Gamma : = \sup_{t} \| \tu_{1} (t, \cdot) \|_{L^2(0, L)} + 1. 
$$

Set 
$$
\tu_2 (t, x) =  W(t, x) \mbox{ in }   (0, + \infty) \times (0, L), 
$$
\begin{equation*}
u_d = \eps \tu_1 + \eps^2 \tu_2 - u \mbox{ in } (0, + \infty) \times (0, L), 
\end{equation*}
\begin{equation*}
f_d = u u_x - \eps^2 \tu_1 \tu_{1, x} \mbox{ in } (0, + \infty) \times (0, L). 
\end{equation*}

We have 
\begin{equation}\left\{
\begin{array}{cl}
u_{d, t} (t, x) + u_{d, x} (t, x) + u_{d, xxx} (t, x)  = f_d (t, x) &  \mbox{ for } t \in (0, +\infty), \, x \in (0, L), \\[6pt]
u_d(t, x=0) = u_d(t, x=L) = u_d(t , x= L)=  0 & \mbox{ for } t \in (0, +\infty), \\[6pt]
u_d(t = 0, \cdot)  = 0 &  \mbox{ in } (0, L). 
\end{array}\right.
\end{equation}

Denote 
\begin{equation}\label{pro-opt-def-g}
g_d = \eps^3 ( \tu_1 \tu_{2, x} +  \tu_2 \tu_{1, x} ) + \eps^4  \tu_2 \tu_{2, x}. 
\end{equation}
We write $f_d$ under the form 
\begin{align*}
f_d =  & (u- \eps \tu_1 - \eps^2 \tu_2) u_x +  (\eps \tu_1 + \eps^2 \tu_2)  (u- \eps \tu_1 - \eps^2 \tu_2)_x   + g_d \\[6pt]
= & - u_d u_x -  (\eps \tu_1 + \eps^2 \tu_2)  u_{d, x}   + g_d. 
\end{align*}
Multiplying the equation of $u_d$ with $u_d$ (which is real), integrating by parts in $(1, t) \times (0, L)$, and using the form of $f_d$ just above  give
\begin{multline}\label{pro-opt-p1}
\int_0^L |u_d(t, x)|^2 \, dx \le \int_0^L |u_d(1, x)|^2 \, dx +  2 \int_1^t  \int_0^L |u_d|^2 |u_x| \, dx \, ds  \\[6pt]
+ \int_1^t \int_{0}^L ( \eps |\tu_1| + \eps^2 |\tu_2| )_x |u_d|^2 \,dx \, ds +  2 \int_1^t \int_0^L |g_d| |u_d|. 
\end{multline}
Since 
$$
\int_0^L |u(t, x)|^2 \, dx  \le \int_0^L |u(0, x)|^2 \le C \eps \mbox{ for } t \ge 0, 
$$
and the effect of the regularity, one has 
\begin{equation}\label{pro-opt-p2}
|u(t,  x)| + |u_x(t,  x)|   \le   C \eps \mbox{ for } t \ge 1, \, x \in [0, L]. 
\end{equation}

Let $a$ be a (small) positive constant defined later (the smallness of $a$ depending only on $L$). Let $t_0 \in [1, a/\eps]$ be such that 
$$
\int_0^L |u_d(t_0, x)|^2 \, dx = \max_{t \in [1, a/\eps]} \int_0^L |u_d(t, x)|^2 \, dx. 
$$
Combining \eqref{pro-opt-p1} with $t = t_0$ and \eqref{pro-opt-p2} yields 
$$
 \int_0^L |u_d(t_0, x)|^2 \, dx \le  \int_0^L |u_d(1, x)|^2 \, dx +  C a  \int_0^L |u_d(t_0, x)|^2 \, dx +    \int_1^{a/\eps} \int_0^L \eps^{-1} |g_d|^2 \, d x. 
$$
This implies, if $a$ is sufficiently small,  
$$
\int_0^L |u_d(t_0, \cdot)|^2 \, dx  \, dx \le C \int_0^L |u_d(1, x)|^2 \, dx  +  C \eps^4
$$
by \eqref{pro-opt-def-g}.

On the other hand, one has 
$$
\int_0^L |u_d(t, \cdot)|^2 \, dx  \, dx \le C \int_0^L |u_d(0, x)|^2 \, dx + C \eps^4  \mbox{ for } t \in [0, 1]. 
$$
We have just proved that, for $a$ sufficiently small, 
$$
\sup_{t \in [0, a/\eps]} \| u_d(t, \cdot) \|_{L^2(0, L)} \le C \Big( \| u_d(0, \cdot) \|_{L^2(0, L)} + \eps^2  \Big). 
$$

Continuing this process, we obtain 
\begin{equation}\label{pro-opt-p3}
\sup_{t \in [0, a n /\eps]} \| u_d(t, \cdot) \|_{L^2(0, L)}  \le  C^n  \| u_d(0, \cdot) \|_{L^2(0, L)} +  \sum_{k=1}^n  C^k \eps^2. 
\end{equation}

We now consider $u$ with 
$$
u(0, \cdot) = \eps \tu_1(0, \cdot) + \eps_2 \tu_2 (0, \cdot). 
$$
Thus 
\begin{equation}\label{pro-opt-p4}
u_d(0, \cdot)  = 0. 
\end{equation}

Fix $ \gamma > 0$  such that 
\begin{equation}
\inf_{t \in \mR} \int_{0}^L |\tu_1(t, x)|^2 \,dx \ge 4 \gamma.  
\end{equation}
With $n $ being the largest integer number such that $C^{n+1} \le \gamma \eps^{-1}$ (we assume now and later on that  $C \ge 2$),  we derive from \eqref{pro-opt-p3} and \eqref{pro-opt-p4} that 
\begin{equation*}
\sup_{t \in [0, a n /\eps]}  \| u_d(t, \cdot) \|_{L^2(0, L)}  \le   \gamma \eps.  
\end{equation*}
Since 
$$
u_d = \eps \tu_1 + \eps^2 \tu_2 - u, 
$$
by the choice of $\gamma$, we have, for $\eps $ sufficiently small,  
\begin{equation*}
 \| u(a n/ \eps, \cdot) \|_{L^2(0, L)}  \ge  \gamma \eps.  
\end{equation*}
We deduce that, with $\tau  = a n/ \eps \sim \eps^{-1} \ln \eps^{-1}$ (hence $\eps^{-1} \sim \tau / \ln \tau$),  
\begin{equation*}
\| u(\tau, \cdot) \|_{L^2(0, L)} \ge \gamma \eps  \ge C \gamma   \ln \tau/ \tau.  
\end{equation*}
The proof is complete. 
\qed

\appendix

\section{Proof of \Cref{lem-Q}}
Let $L \in \cN$ and $z \in \cP_L$. Then, from \cite{Rosier97}, $z = p_m$ for some $1 \le m \le n_L$ and 
$$
\lambda_j = \eta_{j, m}. 
$$
One can then check that $\det Q = 0$. On the other hand, if $z  \neq \pm  2/ (3 \sqrt{3})$  and $\det Q (z) = 0$, it follows that there exists $(a_1, a_2, a_3) \in \mC^3 \setminus \{0 \}$ such that the function $\xi$ defined by 
$$
\xi(x) = \sum_{j=1}^3 a_j e^{\lambda_j (z) x}
$$
satisfies 
$$
\xi(0) = \xi(L) = \xi'(L) = 0. 
$$
Since $\xi''' + \xi' = i z \xi$,  by an integration by parts, one has 
$$
\xi'(0) = 0 \mbox{ if $z$ is real}.  
$$
Hence, from \cite{Rosier97}, if $z \in \mR \setminus  \{  \pm  2/ (3 \sqrt{3})\}$ and $\det Q (z) = 0$, then 
$L \in \cN$ and $z = p_m$ for some $1 \le m \le n_L$. 
We finally
note that, $\{ \pm 2  / (3\sqrt 3) \} \cap \cP_L = \emptyset$ for all $L \in \cN$ 
since, for $k \ge l \ge 1$,
\[
0 \le \frac{(2k + l)(k-l)(2 l + k)}{3 \sqrt{3}(k^2 + kl + l^2)^{3/2}} = \frac{(2k + l)(k^2  + k l - 2 l^2)}{3 \sqrt{3}(k^2 + kl + l^2)^{3/2}} <  \frac{(2k + l)}{3 \sqrt{3}(k^2 + kl + l^2)^{1/2}} < \frac{2}{3 \sqrt{3}}.
\]
The proof is complete. 
\qed

\section{Scilab program for checking $s(k, l) \neq 0$}
\begin{lstlisting}
clc
N=2000;
t=100; 
a=0; 
b=0; 
for k=2:N 
    for l=1:k-1
        h1 = 3 * (k*k + l*l + k*l); 
        h0= 2 * (2*k + l)* (2*l + k)*(k-l); 
        p =  poly([h0 -h1 0 1],'x','c'); 
        r = roots(p); 
        c=exp(4 * %pi * %i * (k-l)/3);
        a1= c* exp(2 * %i * %pi *r(1) /3) 
        			+ exp(-2 * %i * %pi * r(1)/3); 
        a2 = c* exp(2 *  %i * %pi *r(2) /3) 
        			+  exp(-2 *  %i *%pi * r(2)/3);
        a3=c* exp(2 * %i * %pi *r(3) /3) 
        			+ exp(-2 * %i * %pi * r(3)/3);
        s  = r(1)* (r(3) - r(2))* a1 
        			+ r(2)* (r(1) - r(3))* a2 
        			+ r(3)* (r(2) - r(1))* a3;
        if abs(s) < t then t=abs(s); a=k; b=l; 
        end
end
end
disp(a, b, t);
\end{lstlisting}

The outcome  is $t= 0.0000164$, $a= 736$, and $b = 611$. This means 
$$
\min \Big\{|s(k, l)|; 1\le l < k \le 2000 \Big\} = t = 0.0000164
$$ 
and 
$$
s(736,611) = t. 
$$

\medskip 
\noindent \textbf{Acknowledgments.} The author thanks Jean-Michel Coron for interesting discussions and useful comments. 
He also thanks Fondation des Sciences Math\'ematiques de Paris (FSMP) for the Chaire d'excellence which allows him to visit  Laboratoire Jacques Louis Lions  and Mines ParisTech. Part of this work has been done during this visit.

\medskip 
\providecommand{\bysame}{\leavevmode\hbox to3em{\hrulefill}\thinspace}
\providecommand{\MR}{\relax\ifhmode\unskip\space\fi MR }
\providecommand{\MRhref}[2]{%
  \href{http://www.ams.org/mathscinet-getitem?mr=#1}{#2}
}
\providecommand{\href}[2]{#2}

\end{document}